\documentclass{amsart}
\usepackage{amssymb,amsmath,array}

\newtheorem{Lemma}{Lemma}[section]\newcommand{\bel}{\begin{Lemma}}\newcommand{\eel}{\end{Lemma}}
\newtheorem{Proposition}[Lemma]{Proposition}\newcommand{\bprop}{\begin{Proposition}}\newcommand{\eprop}{\end{Proposition}}
\newcommand{\bpr}{{\it Proof. ~}}\newcommand{\epr}{$\blacksquare$\\}
\newtheorem{Property}[Lemma]{Property}\newcommand{\bpropt}{\begin{Property}}\newcommand{\epropt}{\end{Property}}
\newtheorem{Theorem}[Lemma]{Theorem}\newcommand{\bthe}{\begin{Theorem}}\newcommand{\ethe}{\end{Theorem}}
\newtheorem{Remark}[Lemma]{Remark}\newcommand{\beR}{\begin{Remark}\rm}\newcommand{\eeR}{\end{Remark}}
\newtheorem{Definition}[Lemma]{Definition}\newcommand{\bed}{\begin{Definition}}\newcommand{\eed}{\end{Definition}}
\newtheorem{Definition-Proposition}[Lemma]{Definition-Proposition}
\newcommand{\bdp}{\begin{Definition-Proposition}}\newcommand{\edp}{\end{Definition-Proposition}}
\newtheorem{Example}[Lemma]{Example}\newcommand{\bex}{\begin{Example}\rm}\newcommand{\eex}{\end{Example}}
\newtheorem{Corollary}[Lemma]{Corollary}\newcommand{\bcor}{\begin{Corollary}\rm}\newcommand{\ecor}{\end{Corollary}}
\newtheorem{Fact}[Lemma]{Fact}\newcommand{\bfact}{\begin{Fact}\rm}\newcommand{\efact}{\end{Fact}}

\newcommand{\beq}{\begin{equation}}\newcommand{\eeq}{\end{equation}}
\newcommand{\bem}{\begin{displaymath}}\newcommand{\eem}{\end{displaymath}}
\newcommand{\beqa}{\begin{eqnarray}}\newcommand{\eeqa}{\end{eqnarray}}
\newcommand{\bee}{\begin{enumerate}}\newcommand{\eee}{\end{enumerate}}
\newcommand{\bei}{\begin{itemize}}\newcommand{\eei}{\end{itemize}}
\newcommand{\bet}{\begin{tabular}{cccccccc}}\newcommand{\eet}{\end{tabular}}
\newcommand{\bpm}{\begin{pmatrix}}\newcommand{\epm}{\end{pmatrix}}
\newcommand{\bM}{\begin{matrix}}\newcommand{\eM}{\end{matrix}}
\newcommand{\ber}{\begin{array}{l}}\newcommand{\eer}{\end{array}}

\newcommand{\di}{\partial}
\newcommand{\ra}{\!\!\rightarrow\!\!}
\newcommand{\da}{\downarrow}
\newcommand{\tinyM}{\scriptstyle}\newcommand{\tinyT}{\scriptsize}\newcommand{\tinyA}{\tinyM\tinyT}
\newcommand{\smallA}{\footnotesize}

\newcommand{\cO}{{\mathcal{O}}}

\newcommand{\mC}{\mathbb{C}}

\newcommand{\mP}{\mathbb{P}}\newcommand{\tP}{\check{\mP}}\newcommand{\mS}{\mathbb{S}}\newcommand{\mXL}{\mP_x^2\times\tP^2_l}
\newcommand{\mR}{\mathbb{R}}
\newcommand{\mZ}{\mathbb{Z}}

\newcommand{\al}{\alpha}\newcommand{\be}{\beta}\newcommand{\De}{\Delta}\newcommand{\de}{\delta}
\newcommand{\ga}{\gamma}\newcommand{\Ga}{\Gamma}
\newcommand{\ep}{\epsilon}\newcommand{\si}{\sigma}
\newcommand{\Si}{\Sigma}

\newcommand{\tSi}{\widetilde{\Sigma}}
\newcommand{\tv}{\tilde{v}}

\newcommand{\tAux}{\widetilde{\rm Aux}}

\newcommand{\bSi}{\bar\Sigma}
\newcommand{\lSi}{{\overline{\Sigma}}}\newcommand{\ltSi}{{\overline{\widetilde\Sigma}}}

\newcommand{\ol}[1]{\overline{#1}}
\newcommand{\li}{~\\ $\bullet$ }\newcommand{\ls}{~\\ $\star$ }
\newcommand{\into}{\stackrel{i}{\hookrightarrow}}

\newcommand{\bin}[2]{{#1\choose{#2}}}

\newcommand{\gNnd}{generalized Newton-non-degenerate }\newcommand{\Nnd}{Newton-non-degenerate}

\newcommand{\ND}{Newton diagram}\newcommand{\D}{diagram }

\newcommand{\sset}{\subset}\newcommand{\smin}{\setminus}
\newcommand{\sx}{{\mS_x}}\newcommand{\sy}{{\mS_y}} \newcommand{\lxy}{{\overline{xy}}}

\newcommand{\omp}{ordinary multiple point}
\newcommand{\wrt}{with respect to }\newcommand{\mult}{{\rm mult} }\newcommand{\jet}{{\rm jet} }
\renewcommand{\kill}{{\rm kill} }
\newcommand{\Aux}{{\rm Aux} }

\newcommand{\suml}{\sum\limits}\newcommand{\capl}{\mathop\cap\limits}

\newcommand{\bl}{\langle}\newcommand{\br}{\rangle}

\newcommand{\mesh}[7]
{
\put(#1,#2){\vector(1,0){#6}}  \put(#1,#2){\vector(0,1){#7}}
\setcounter{tempx}{#3}  \addtocounter{tempx}{1}  \setcounter{tempy}{#4}  \addtocounter{tempy}{1}
\multiput(#1,#2)(#5,0){\value{tempx}}{\multiput(-1.5,-0.5)(0,#5){\value{tempy}}{.}}
}

\title{O\MakeLowercase{n the enumeration of complex plane curves with two singular points}}
\author[]{D\MakeLowercase{mitry} K\MakeLowercase{erner}\\\small B\MakeLowercase{en} G\MakeLowercase{urion}
U\MakeLowercase{niversity}, I\MakeLowercase{srael}}\date{\today}
\address{Department of Mathematics, Ben Gurion University of the Negev, P.O.B. 653, Be'er Sheva 84105, Israel.}
\email{kernerdm@math.bgu.ac.il}

\thanks{\date{\today}This paper is an updated version of the preprint: Max Planck Institut f\"ur Mathematik Bonn,
Germany, MPIM2007-76.
\\
The author was supported by the Max Planck Institut f\"ur Mathematik (Bonn, Germany) and the Skirball postdoctoral fellowship of the Center
of Advanced Studies in Mathematics at the Mathematics Department of Ben-Gurion University (Be'er-Sheva, Israel).}



\begin{document}\setcounter{secnumdepth}{6} \setcounter{tocdepth}{1}\newcounter{tempx}\newcounter{tempy}
\begin{abstract}
We study equi-singular strata of plane curves with two singular points of prescribed types.
The method of the  previous work \cite{Ker06 enum.unising.curves} is generalized to this case.
In particular we consider the enumerative problem for
 plane curves with two singular points of {\it linear} singularity types.

First the problem for two \omp s of fixed multiplicities is solved. Then the enumeration for arbitrary
linear types is reduced to the case of \omp s and to the understanding of "merging" of singular points.
Many examples and numerical answers are given.
\end{abstract}

\maketitle

\tableofcontents
\section{Introduction}
\subsection{The setup and the problem}
We work with (complex) algebraic curves in $\mP^2$. Consider the parameter space of degree-$d$ curves,
i.e. the complete linear system $|\cO_{\mP^2}(d)|$. It is a projective space of dimension $N_d=\bin{d+2}{2}-1$.

A classical enumerative problem is:
{\it given the singularity types $\mS_1\dots\mS_r$, "how many" curves of the linear system $|\cO_{\mP^2}(d)|$ possess
singular
points of these types?}  (To make this number finite one imposes a sufficient amount of generic base points. The
degree $d$ is assumed big enough to avoid various pathologies. By the singularity type in this paper we
always mean the local, embedded, topological singularity type
of a reduced plane curve. For more details and related notions cf.\S\ref{Sec Background}.)

In this paper we consider the case of two prescribed singular points (the case of one singular point
was considered in \cite{Ker06 enum.unising.curves}). First reformulate the problem.
The parameter space $|\cO_{\mP^2}(d)|$ is stratified according to the  singularity types of curves.
The generic point of $|\cO_{\mP^2}(d)|$
corresponds to an (irreducible, reduced) smooth curve. The set of points corresponding to the singular curves
is called the {\it discriminant} of plane curves. It is an irreducible projective hypersurface in $|\cO_{\mP^2}(d)|$.

\bed
For the given  singularity types $\mS_1\dots \mS_r$, the {\it equisingular stratum}
$\Si_{\mS_1\dots \mS_r}\subset|\cO_{\mP^2}(d)|$ is the set of points corresponding to the curves with $\mS_1\dots \mS_r$
singularities.
\eed
 The generic point of the discriminant lies in the stratum of nodal curves
($\Si=\lSi_{A_1}$). Other strata correspond to higher singularities. The stratum of $\de$-nodal curves $\Si_{\de A_1}$
 (whose closure $\lSi_{\de A_1}$ contains the stratum of curves of a given genus) is the classical Severi variety.
Other strata are $\Si_{A_k}$, $\Si_{D_k}$, $\Si_{E_k}$ etc.
The strata are quasi-projective varieties \cite{Greuel-Lossen96}.
For a comprehensive introduction to these
equi-singular families and related notions cf. \cite{GLSbook1},\cite{GLSbook2}.

For small $d$ various pathologies can occur,
but for sufficiently high degrees (given the singularity types $\mS_1\dots \mS_r$) the strata
are non-empty, irreducible, smooth, of expected (co-)dimension.
One sufficient condition for this is \cite[\S I.3]{Dim}: $d\ge\sum o.d.(\mS_i)+r-1$, here $o.d.$ are
the orders of determinacy. There are also other sufficient conditions, formulated in terms of the Milnor or Tjurina numbers of
the types or the $\ga$ invariant (cf.  \cite[\S 3]{GLS06} and \cite[\S IV.2]{GLSbook2}), e.g.
$\sum \ga(\mS_i)\le(d+3)^2$. In this paper we always assume $d$ big enough.
\\\\
By definition each stratum is embedded into $|\cO_{\mP^2}(d)|$, thus its natural compactification is just the
topological closure.
The closures of the strata are singular (often in co-dimension 1).
The closed stratum defines the homology class $[\lSi_{\mS_1\dots \mS_r}]\in H_*(|\cO_{\mP^2}(d)|,\mZ)\approx\mZ$
 in the homology of the corresponding dimension.
The degree of this class is the degree of the stratum $\deg(\lSi_{\mS_1\dots \mS_r})$,
obtained by the intersection with the generic plane in $|\cO_{\mP^2}(d)|$ of the complementary dimension.
This corresponds to imposing generic base points.
So this degree is {\it "the number" of curves possessing the prescribed singularities}.
\\\\
In \cite{Ker06 enum.unising.curves} we proposed a method to compute the degrees of strata $\lSi_{\mS}$ for curves
with just one singular point.
In \cite{Ker05Hypersurfaces} the method was generalized to some singular hypersurfaces in $\mP^n$.
The current paper is an application of the method to the enumeration of curves with two singular points,
i.e. to compute the classes $[\lSi_{\mS_1\mS_2}]$.
\\\\
{\bf Acknowledgements.}
This work is a tail of my PhD, done under the supervision of E.Shustin, to whom I wish to express my deepest gratitude.
The conversations and advices of P.Aluffi, G.-M.Greuel, I.Tyomkin were highly important.
Careful reading of two anonymous referees has improved the text significantly.

Most of the work was done during my stay in Max Planck Institut f\"ur Mathematik, Bonn in 2006-2007. Many thanks
for excellent working conditions.
\subsection{The known  results}\label{SecTheResultsOldAndNew}
Since the question is completely classical, lots of {\it particular} results are known. We mention only a few
(for a much better discussion cf. \cite{Klei76,Klei85}, \cite[Remark 3.7]{KleiPien01} and \cite{Kaz01}).
\li The most classical case is $\deg\lSi_{\de A_1}$. The history of solved cases starts from $\de=1$ in \cite{Steiner1848},
$\de=2$ in \cite{Cayley1866} and $\de=3$ in \cite{Roberts1875}. The few-nodal cases were recalculated many times
by various methods, e.g. for $\de\le3$ in \cite{HarrPand95}, for $\de\le 5$ in \cite{Vain81,Vain95}, for $\de\le 8$ in
\cite{KleiPien98,KleiPien01,KleiPien04}.
Various general algorithms have been given (e.g. \cite{Ran89,Ran02},\cite{CapHar98-1,CapHar98-2}).
It was conjectured \cite{Goett98} that $\deg\lSi_{\de A_1}$ is a polynomial in $d$ of degree $2\de$, for $\de<2d-1$.
This was proved in \cite{Liu2000}, where a general method of calculation of $\deg\lSi_{\de A_1}$ was also proposed.
Another proof has appeared recently: \cite{Fomin-Mikhalkin09}.

Finally, some many-nodal cases were treated, as they amount to enumeration of low genus plane curves
(with Gromov-Witten invariants, quantum cohomology, \cite{KontMan94} and all that.)
The later approach is effective for low genus curves, i.e. when the
number of nodes (or higher singularities) is almost maximal ($\de\lesssim\frac{(d-1)(d-2)}{2}$).
It seems to be non-effective for high genus computations
(e.g. in the case of just a few singular points).
\li The $\deg\lSi_{A_2}$ was predicted in \cite{Enriques32} and (re)proved
in \cite[pg.151]{Lascoux77}, \cite[8.6.3, pg.416]{Vain81}, \cite[pg.86-88]{diFrItz95} and \cite[pg.10]{Alufi98}. In the
later paper $\deg\lSi_{A_3}$ was also obtained.
In \cite[theorem 1.2]{KleiPien98} the degrees of $\lSi_{D_4,A_1}$, $\lSi_{D_4,2A_1}$, $\lSi_{D_4,3A_1}$,
$\lSi_{D_6,A_1}$ and $\lSi_{E_7}$ were computed.
The $\deg(\lSi_{A_2A_1})$ was computed in \cite[8.6.4, pg.416]{Vain81}
\li The impressive breakthrough has been recently achieved in \cite{Kaz01}-\cite{Kaz03-2}. The proposed topological method
allows (in principle) to compute the degree of {\it any} stratum (with lots of explicit results in \cite{Kaz03-hab}).
In particular he presents the answers $\deg\lSi_{\mS_1\dots \mS_r}$ for the total codimension: $\sum codim(\mS_i)\le7$.
For higher singularities the method is not quite efficient, as it solves the problem {\it simultaneously} for all the
singularity types of a given co-dimension. So, first one should classify the singularities (by
now the classification seems to exist up to codimension 16 only). Even if this is done, one faces the problem
of enumerating huge amount of cases (the number of types grows exponentially with the codimension).
And of course, each computation can give a result for a specific choice of singularity types,
it is not clear whether the method allows to obtain results for some {\it series} of singularities.
\li In \cite{Ker06 enum.unising.curves} the problem was solved for curves with one singular point of
 an {\it arbitrary} given singularity type. The proposed method gives immediate answer (explicit formulas)
 for some specific series of types
(the so-called {\it linear}).
For all other (series of) types it gives an explicit algorithm.
\\
\\
\\
Despite numerous isolated results as above it seems that currently there is no  effective universal
method to compute $\deg\bSi_{\mS_1\dots \mS_k}$. Not speaking about a general
formula giving  this degree for various choices.
And most methods use the information about adjacency of types: which singularity is obtained when the
singular points of the types $\mS_1\dots \mS_k$ are merged (collide) generically?
\subsection{Our results}\label{Sec Intro Our Results}
We consider the equisingular strata of curves with two singular points: $\lSi_{\sx\sy}$. We restrict mostly to the case
of linear singularity types (cf. definition \ref{Sec Background Linear Singularities}). The simplest examples
of linear singularity types are $x_1^p+x_2^{q}$ for $p\le q\le 2p$.
(In particular they include \omp, $A_{k\le3}$, $D_{k\le6}$, $E_{k\le8}$,
$J_{10}$, $Z_{k\le13}$ etc.)
\\
\\
\parbox{13.2cm}
{We construct the partial resolutions $\ltSi_{\sx\sy}\to\lSi_{\sx\sy}$ as on the diagram.
Here $\Aux$ is an auxiliary space tracing parameters of the singularities,
e.g. the singular points, the lines of the tangent cones, the line $l=\overline{xy}$.
The cohomology class $[\lSi_{\sx\sy}]\in H^*(|\cO_{\mP^2}(d)|)$ is completely determined by the class
$[\ltSi_{\sx\sy}]\in H^*(|\cO_{\mP^2}(d)|\times \Aux)$, so we compute the later.}\hspace{0.6cm}
$\bM \ltSi_{\sx\sy}\subset&|\cO_{\mP^2}(d)|\times \Aux
\\\da&\da
\\\lSi_{\sx\sy}\subset&|\cO_{\mP^2}(d)|
\eM$
\\
In \S\ref{Sec Example Two OMP's} the enumerative problem for two \omp s is solved, i.e.
the cohomology class of $\ltSi_{x^{p+1}_1+x^{p+1}_2,y^{q+1}_1+y^{q+1}_2}$ is computed. The enumeration is
done in two ways: stepwise intersection with hypersurfaces and degeneration to reducible curves.

In \S\ref{Sec Method Degenerations To OMP}
we use the method of degenerations from
\cite{Ker06 enum.unising.curves} to express the class $[\ltSi_{\sx\sy}]$ (for $\sx,\sy$ arbitrary linear singularities)
via
$[\ltSi_{x^{p+1}_1+x^{p+1}_2,y^{q+1}_1+y^{q+1}_2}]$ and some necessary classes $[\ltSi_\mS]$
(the classes of some strata of curves with one singular point of type $\mS$).
\\
\\
The strata $\ltSi_\mS$ appear inevitably, they are irreducible components of the restriction to the diagonal:
$\ltSi_{\sx\sy}|_{x=y}$. Therefore the enumerative problem is reduced to the description of
such "collision of singular points". This later problem is very complicated in general.
But it is manageable when the types $\sx,\sy$ are {\it linear}. In \cite{Ker07Collisions}
we gave an explicit algorithm to classify the "results of collisions", i.e. the irreducible components of
$\ltSi_{\sx\sy}|_{x=y}$. In fact the algorithm even gives the defining ideals of the corresponding
strata, this is briefly recalled in \S\ref{Sec Background Collisions}. So, each class $[\ltSi_\mS]$ is
computable  by methods of \cite{Ker06 enum.unising.curves}.
\\
\\
\\
So, for linear singularities, the method allows to compute the class $[\ltSi_{\sx\sy}|_{x=y}]$.
Many numerical answers are given in the Appendix.
As mentioned above, the results are valid for $d\gg0$. In \S\ref{Sec Method Validitiy of Results}
we give a cheap sufficient bound (in terms of the orders of determinacy of $\sx\sy$), though the bound
is far from being necessary.

For small $d$ (but big enough so that the strata $\ltSi_{\sx\sy},\ltSi_{\mS}$ are irreducible, of expected dimension)
the problem is still well defined, but our formulas are inapplicable.
For each such $d$ the enumeration should be done separately, cf. \cite[\S5]{Ker06 enum.unising.curves}.
\\
\\
\\
Our approach is most naive and classical. In some
sense it is a brute-force calculation. Correspondingly it is often long and cumbersome. The advantages
of the method are:
\li The method gives a recursive algorithm, consisting of routine parts.
\li The computation does not assume any preliminary classification of singularities.
\li The method seems to be more effective than other approaches (to the best of our knowledge).
In particular, in Appendix we present the results for some {\it series} of types (as compared
to {\it single, isolated} results previously known).
Thus we can treat the question:
 how do the degrees of the strata depend on the invariants of the singularity types for some series of singularities?

In the case of one singular point many examples suggested in \cite{Ker06 enum.unising.curves} that (at least for linear types) the
degrees depend {\it algebraically} on the parameters of series (e.g. for \omp~there is a polynomial
dependence on the multiplicity).

One might conjecture that in the multi-singular case the dependence will be also algebraic.
Unfortunately, the simplest case already provides a counterexample: for the two \omp s of multiplicities $p+1,q+1$
the degree of the corresponding stratum depends on $max(p,q)$ and $min(p,q)$, cf. remark
\ref{Rem Class depends non algebraically}.
\li As the final result we obtain the multi-degree of the (partial) resolution of the stratum $\lSi_{\sx\sy}$.
(The actual degree is just a particular coefficient in a big polynomial.) This multi-degree contains
many important numerical invariants, e.g. enumeration of ($\sx,\sy$) with one or two singular points restricted
to some curves, or with some conditions on the tangents to the branches. More generally: when
the parameters of the singularity (the points, the tangents)
are restricted to a subvariety of the original parameter space. So, this solves a whole class
of related enumerative problems.
\subsection{Contents of the paper}
In \S\ref{Sec Background}  we fix the notations and recall some necessary notions. In particular,
{\it linear singularities}
are introduced in \S \ref{Sec Background Linear Singularities} and {\it collisions} of singular points
are discussed in \S \ref{Sec Background Collisions}.
In \S \ref{Sec Method} we describe the method. First we recall
in \S\ref{Sec Method Strata of Unisingular Curves} the case of curves with one singular point.
In \S\ref{Sec Method Strata For two singular points} we discuss the case of two singular points.
\S\ref{Sec Method Stepwise Intersection} presents the direct approach: obtaining the cohomology class
as the chain of intersections. \S\ref{Sec Method Degenerations To OMP} presents the indirect
approach: degeneration to \omp s of higher multiplicities.

The method is applied to several particular cases in \S\ref{Sec Examples},
illustrating the computation.
Explicit formulas are given in Appendix, for the details of numerical computation see the attached
Mathematica file \cite{Ker07Bisingular}.
\section{Some relevant notions and auxiliary results}\label{Sec Background}
\subsection{Notations and some classical facts}
\subsubsection{Coordinates and variables}\label{Sec Background Coordinates}
We work with various projective spaces and their subvarieties. Adopt the following notation.
\\
A curve is denoted by $C\sset\mP^2$ or by the defining polynomial $f$, the parameter space of such
curves of degree $d$ is $|\cO_{\mP^2}(d)|$.
Let $x\in\mP^2_x$ be a point in the projective plane, the homogeneous coordinates are $(x_0,x_1,x_2)$.
The generator of the cohomology
ring of this $\mP^2_x$ is denoted by the upper-case letter $X$, so
that $H^*(\mP^2_x,\mZ)=\mZ[X]/(X^{3})$.
By the same letter we denote the class of a line in the homology of $\mP^2_x$. Since it is always clear,
where we speak about coordinates and where about (co)homology classes, no confusion arises.
For example, consider the hypersurface
\beq
V=\{(x,y,f)|~f(x,y)=0\}\subset\mP^2_{x}\times\mP^2_{y}\times|\cO_{\mP^2}(d)|
\eeq
Here $f$ is a bi-homogeneous polynomial of bi-degree $d_x,d_y$ in homogeneous coordinates
$(x_0,x_1,x_2)\in\mP^2_x$, \mbox{$(y_0,y_1,y_2)\in\mP^2_y$}.
The coefficients of $f$ are the homogeneous coordinates on the parameter space $|\cO_{\mP^2}(d)|$.
The cohomology class of this hypersurface is
\beq\label{EqDemostratCohomClasses}
[V]=d_xX+d_yY+F\in H^2(\mP^2_{x}\times\mP^2_{y}\times|\cO_{\mP^2}(d)|,\mZ)
\eeq
A (projective) line through the point $x\in\mP^2_x$ is defined by a 1-form $l$ (so that $l\in\tP^2_l,~l(x)=0)$.
Correspondingly the generator of $H^*(\tP^2_l)$ is denoted by $L$.
\\\\
For a variety $V\into |\cO_{\mP^2}(d)|\times\mP^2_x\times\mP^2_y$ the pullbacks $i^*(X),i^*(Y),i^*(F)$ are constantly used.
To simplify the formulas we
denote the pulled back classes by the same letters ($X,Y,F$).
This brings no confusion as the ambient space of intersection is always specified.
\\\\
We often work with symmetric $p-$forms $\Omega^{(p)}\!\!\in\!\! Sym^p(\check{V}_3)$, here $\check{V}_3$ is
a 3-dimensional vector space of linear forms.
Fix coordinates in $\check{V}_3$, then the form is a symmetric tensor with $p$ indices
($\Omega^{(p)}_{i_1,\dots,i_p}$). We write
$\Omega^{(p)}(\underbrace{x\dots x}_{k})$ as a shorthand for the tensor, multiplied $k$ times by the point $x\in V_3$:
\beq
\Omega^{(p)}(\underbrace{x\dots x}_{k}):=\sum_{0\le i_1,\dots,i_k\le2}\Omega^{(p)}_{i_1,\dots,i_p}x_{i_1}\dots x_{i_k}
\eeq
So, for example, the expression $\Omega^{(p)}(x)$ is a $(p-1)-$form. Unless stated otherwise, we assume the
symmetric form $\Omega^{(p)}$ to be generic (in particular non-degenerate, i.e. the corresponding curve
$\{\Omega^{(p)}(\underbrace{x\dots x}_{p})=0\}\subset\mP^2_x$ is smooth).

Symmetric forms occur typically as tensors of derivatives of order $p$ in homogeneous coordinates: $f^{(p)}$.
Sometimes, to emphasize
the point at which the derivatives are calculated we assign it. So, e.g.
$f|_x^{(p)}(\underbrace{y,\dots,y}_k)$ means a symmetric $(p-k)$ form: the tensor of derivatives of
order $p$, calculated at the point $x$, and multiplied $k$ times with $y$.

Sometimes we address a particular component of the tensor, e.g.
$(f|_x^{(p)})_{i_1\dots i_p}:=$ $\frac{\di}{\di x_{i_1}}\dots \frac{\di}{\di x_{i_p}}f|_x\equiv$ $\di_{i_1}\dots .\di_{i_p}f|_x$.

The Euler identity $f|_x^{(1)}(x)\equiv\sum x_i\di_if|_x=df$ and its consequences
(e.g. $f|_x^{(p)}\underbrace{(x\dots x)}_{k}\sim f|_x^{(p-k)})$ are tacitly assumed.
\subsubsection{Blowup along the diagonal.}\label{Sec Background Blowup Diagonal}
The diagonal $\Delta=\{x=y\}\subset\mP_x^n\times\mP_y^n$ appears constantly.
Its class is \cite[example 8.4.2]{Ful98}
\beq
[\Delta]=\sum^n_{i=0} X^{n-i}Y^i\in H^{2n}(\mP_x^n\times\mP_y^n,\mZ)
\eeq
For example, the proportionality condition of two symmetric forms $f^{(p)}|_x\sim g^{(p)}|_x$, for $x\in\mP^2$,
 is just the coincidence
of the corresponding points in a big projective space, thus its class is given by the above type formula
(with $n=\bin{p+2}{2}-1$).
\\
\\
The blowup of $\mP_x^2\times\mP_y^2$ over the diagonal $\{x=y\}$ is
described nicely as the incidence variety of triples: a line and
a pair of its points.
\beq
Bl_\De(\mP^2_x\times\mP^2_y)=\{(x,y,l)|~x\in l\ni y\}\into\mP^2_x\times\mP^2_y\times\tP^2_l,~~~E_\De=\{x=y\in l\}\subset Bl_\De(\mP^2_x\times\mP^2_y)
\eeq
\bel
$\bullet$ $Bl_\De(\mP^2_x\times\mP^2_y)\subset\mP^2_x\times\mP^2_y\times\tP^2_l$ is a complete intersection,
with cohomology class: $[Bl_\De(\mP^2_x\times\mP^2_y)]=(L+X)(L+Y)\in H^4(\mP^2_x\times\mP^2_y\times\tP^2_l)$.
\li The class of the exceptional divisor is:
$[E_\De]=(L+X)(X^2+XY+Y^2)\in H^6(\mP^2_x\times\mP^2_y\times\tP^2_l)$ and
$[E_\De]=X+Y-L\in H^2(Bl_\De(\mP^2_x\times\mP^2_y))$.
\eel
The two formulas for $E_\De$ are related by the pushforward $i_*$.
The identity $i_*(X+Y-L)=(L+X)(X^2+XY+Y^2)\in H^6(\mP^2_x\times\mP^2_y\times\tP^2_l)$ is directly verified.
\\
\bpr
Note that $E_\De$ is a transversal intersection of the two conditions ($x=y$ and $l(x)=0$). Therefore:
$[E_\De]=(L+X)(X^2+XY+Y^2)\in H^6(\mP^2_x\times\mP^2_y\times\tP^2_l)$.

To obtain the class of the exceptional divisor in the ring $H^*(Bl_\De(\mP^2_x\times\mP^2_y))$ note that
 the hypersurface $\bpm x_0&x_1\\ y_0&y_1\epm=0$ contains the exceptional divisor $E_\De$ and also
the residual divisor $l_2=0$.
\epr
\subsection{Some notions from singularities}\label{Sec Background SingularityTypes}
\bed\cite[pg.202]{GLSbook1}
Let $(C_x,x)\subset(\mC^2_x,x)$ and $(C_y,y)\subset(\mC^2_y,y)$ be two germs of isolated curve singularities.
They are (locally, embedded, topologically) equivalent if there exist a local homeomorphism
 $(\mC^2_x,x)\mapsto(\mC^2_y,y)$ mapping $(C_x,x)$ to $(C_y,y)$.

 The corresponding equivalence class is
 called the (local embedded topological) singularity type. The variety of points (in the parameter space
$|\cO_{\mP^2}(d)|$), corresponding to the curves with prescribed singularity types $\sx\dots \mS_z$ is called
the { equisingular stratum} $\Si_{\sx\dots \mS_z}$.
\eed
The singularity type can be specified by a polynomial representative of the type.
 Several simplest types have the following representatives
(all the notations are from \cite[\S I.2]{AGLV}, we omit the moduli of analytic classification):
\beq\scriptstyle\ber
A_k:x^2_2+x^{k+1}_1,~~D_k:x^2_2x_1+x^{k-1}_1,~~E_{6k}:x^3_2+x^{3k+1}_1,~~E_{6k+1}:x^3_2+x_2x^{2k+1}_1,~~
E_{6k+2}:x^3_2+x^{3k+2}_1\\
J_{k\ge1,i\ge0}:x^3_2+x^2_2x^k_1+x^{3k+i}_1,~~Z_{6k-1}:x^3_2x_1+x^{3k-1}_1,~~
Z_{6k}:x^3_2x_1+x_2x^{2k}_1,~~Z_{6k+1}:x^3_2x_1+x^{3k}_1\\
X_{k\ge1,i\ge0}:x^4_2+x^3_2x^k_1+x^2_2x^{2k}_1+x^{4k+i}_1,~~W_{12k}:x^4_2+x^{4k+1}_1,~~W_{12k+1}:x^4_2+x_2x^{3k+1}_1
\eer\eeq
For a curve defined by $f=f_p+f_{p+1}+\dots \in\mC\{x_1,x_2\}$ the projectivized tangent cone is the (non-reduced)
scheme $\{f_p=0\}\sset\mP^1$. We denote it as $\mP T_{(C,0)}=l^{p_1}_1\dots l^{p_k}_k$, where $\{l_\al\}$ are the
distinct lines of the tangent cone and $\{p_\al\}$ are their multiplicities, i.e. the local multiplicities of
 the scheme $\{f_p=0\}$. Thus in particular $\sum_\al p_\al=\mult(C)$).

Associated to the tangent cone is the tangential decomposition: $(C,0)=\cup (C_\al,0)$.
 Here the tangent cone of each $(C_\al,0)$ is just one line
and $T_{(C_\al,0)}\neq T_{(C_\be,0)}$, but the germs $(C_\al,0)$ can be further reducible.

\subsubsection{Newton-non-degenerate singularities}
Given the representative $f=\sum a_{\bf I}{\bf x}^{\bf I}\in\mC\{\bf x\}$ of the singularity type, one can draw
the {\it Newton diagram} $\Ga_f$ of the singularity. Namely, one marks the points ${\bf I}$ corresponding to non-vanishing
monomials in $f$, and takes the convex hull of the sets ${\bf I}+\mR_+^2$. The envelope of the convex
hull (the chain of segments) is the Newton diagram.
\bed\label{DefNewtonNonDegenerate}
\li The function is called \Nnd ~ \wrt its diagram if the truncation $f_\si$
 to every face of the diagram is non-degenerate (i.e. the truncated function has no
singular points in the torus $(\mC^*)^2$).
\li The curve-germ is called \gNnd if it can be brought to a Newton-non-degenerate form by a
locally analytic transformation.
\li The singular type is called \Nnd~ if it has a Newton-non-degenerate representative.
\eed
The \Nnd~ type is completely specified by the Newton diagram of (any of) its \Nnd~ representative \cite{Oka79}.

In the tangent cone of the singularity $T_C=\{l_1^{p_1}\dots .l_k^{p_k}\}$, consider the lines appearing with the
multiplicity 1. They correspond to smooth branches, not tangent to any other branch of the singularity.
\bed\label{DefNonFreeBranches}
The branches as above are called {\it free} branches. The tangents to the non-free branches are called {\it non-free tangents}.
\eed

Newton-non-degeneracy implies strong restrictions on the tangent cone: there are at most two non-free tangents.
\bpropt\label{ThmNewtonNonDegenerateTangentCone} Let $T_C=\{l_1^{p_1},\dots .,l_k^{p_k}\}$ be the tangent cone
 of the germ $C=\cup C_j$. If the germ is \gNnd then $p_\al>1$ for at most two tangents $l_\al$.
\epropt
Indeed, suppose $(C,0)$ is \Nnd~in some coordinates, let $f=f_p+f_{p+1}+\dots .$ be its locally defining function.
The non-degeneracy means that $f_p$ has no singular points in $\mC^*$, hence the statement.
\subsubsection{Linear singularities}\label{Sec Background Linear Singularities}
By definition the equisingular stratum $\lSi_\mS$ is a subvariety of $|\cO_{\mP^2}(d)|$. Let the tangent cone be
$T_{(C,0)}=l_1^{p_1}\dots l_k^{p_k}$ where $l_1\dots l_k$ are (distinct) lines passing through $x$.
Consider the subvariety
\beq
\lSi_\mS\supset\lSi_\mS|_{x,\{l_i\}}:=\overline{\{\text{curves of degree $d$ with $\mS$ at $x\in\mP^2$
and $T_{(C,x)}=l_1^{p_1}\dots l_k^{p_k}$}\}}
\eeq
\bed
The singularity (the germ, the stratum, the type) is called linear if $\lSi_\mS|_{x,\{l_i\}}$ is a
linear subspace of $|\cO_{\mP^2}(d)|$ for some (and hence for any) choice of $x,\{l_i\}$.
\eed
For example, an \omp~is linear. For a linear singularity type $\mS$ an appropriate modification
of the equisingular stratum $\ltSi_\mS\to\lSi_\mS$ (defined later)
fibres over the auxiliary space of lines through the points $\Aux=\{(x,l_1\dots l_k)|~\forall i:~x\in l_i\}$.
\\
\\
Linear types are abundant by the following observation. Let $(C,0)=\cup_\al(C_\al,0)$ be the
{\it tangential} decomposition as above. Let $\mS_\al$ be the singularity type of $(C_\al,0)$,
note that the types $\mS_1\dots \mS_k$  are completely determined by $\mS$.
For each non-smooth $(C_\al,0)$ choose the coordinates: one axis is tangent to the branch,
the other axis is generic.
Assume  the \ND~$\Ga_{(C_\al,0)}$ is commode, i.e. intersects all the coordinate axes, and $(C_\al,0)$ is non-degenerate in the chosen coordinates.
For each face $\si$ of the \ND~ let $a_\si$ be the angle between the face and the coordinate axis $\hat{x}_1$.
\bprop\cite[section 3.1]{Ker06 enum.unising.curves}\label{Thm Slope On ND for linear types}
Under the assumptions as above the type $\mS$ is linear iff each $\mS_\al$ is linear.
 And $\mS_\al$ is linear iff every face $\si$ of the \ND~$\Ga_{(C_\al,0)}$ has a
 bounded slope: $\frac{1}{2}\leq|{tg}(a_\si)|\leq2$.
\eprop
\bex
The simplest class of examples of linear singularities is defined by the series: $f=x^p+y^q,~~p\leq q\leq2p$.
In general, for a given series only for a few types of singularities the strata can be linear.
In the low modality cases the linear types are:
\li{Simple singularities (no moduli)}: $A_{1\le k\le3},~~D_{4\le k\le6},~~E_{6\le k\le8}$
\li{Unimodal singularities}: $X_9(=X_{1,0}),~~J_{10}(=J_{2,0}),~~Z_{11\le k\le13},~~W_{12\le k\le 13}$
\li{Bimodal}: $Z_{1,0},~~W_{1,0},~~W_{1,1},~~W_{17},~~W_{18}$
\eex
Most singularity types are nonlinear. For example, if a curve has an $A_4$ point,
then the choice of one axis as the generic tangent brings it to the Newton diagram of $A_3$, but with
the defining function: $(\al x_2+\be x^2_1)^2+\ga x^5_1+\dots .$.
 The stratum of curves whose local defining
 equation begins with expansion of this type (i.e. $\al,\be$ are not fixed) is not a linear subspace in $|\cO_{\mP^2}(d)|$.

\subsubsection{Finite determinacy}\label{Sec Background Finite Determinacy}
The finite determinacy theorem of Tougeron (cf. \cite[\S I.1.5]{AGLV} or \cite[\S I.2.2]{GLSbook1})
states that the topological type of
the  curve-germ is fixed by a finite jet of the defining series. Namely, for every
type $\mS$, there exists $k$ such that for all bigger $n\ge k$: ($\{\jet_n(f)=0\}$ has type $\mS$) implies
($\{f=0\}$ has type $\mS$).
The minimal such $k$
is called: {\it the order of determinacy} (for contact equivalence). E.g. $o.d.(A_k)=k+1$, $o.d.(D_k)=k-1$.
The classical theorem \cite[thm I.2.23]{GLSbook1} reads: if $m^{k+1}\subset m^2Jac(f)$
then $o.d.(f)\le k$.
\subsection{Equisingular strata and related questions}
\subsubsection{Resolution of the singularities of closured strata}\label{Sec Background Improving Fibration}
~\\\parbox{14.5cm}
{The following situation occurs frequently. Let $\Aux$ be a smooth, irreducible projective variety and
 the projection $\tSi\stackrel{\pi}{\to}\Aux$ a locally trivial fibration over the base $\pi(\tSi)\sset \Aux$,
 such that $\ol{\pi(\tSi)}=\Aux$ and the fibres are $\mP^n$
linearly embedded into $|\cO_{\mP^2}(d)|$. In particular $\tSi$ is irreducible and smooth.
}~~~~
$\bM
\tSi\subset\ltSi\sset \Aux\times|\cO_{\mP^2}(d)|\\\da\hspace{0.8cm}\searrow\hspace{0.8cm}\swarrow\pi~~~~~~~~~\\\pi(\tSi)\sset \Aux~~~~~~~~~~~~~~~
\eM$
\\\\
We want to compactify $\tSi$ in a smooth way, preserving the bundle structure. Start from the topological closure
$\ltSi\subset \Aux\times|\cO_{\mP^2}(d)|$. Note that the projection $\ltSi\to \Aux$ is well defined
(being the restriction of $\pi$) and is surjective.
\bprop\label{Thm Desingularization fo Strata}
There exists a birational morphism of smooth varieties $\tAux\stackrel{\phi}{\to}\Aux$ with the properties:
\li It is an isomorphism over $\pi(\tSi)$ and a finite collection of blowups with smooth centers in $\Aux\smin\pi(\tSi)$
\li The corresponding lifting $\ltSi'\sset\tAux\times|\cO_{\mP^2}(d)|$ is a projective fibration over $\tAux$,
with fibres $\mP^n$ linearly embedded into $|\cO_{\mP^2}(d)|$.
\eprop
In particular, $\ltSi'$ is smooth.
\\
\bpr
As the fibres of the projection $\tSi\to\pi(\tSi)\sset \Aux$ are linear subspaces of $|\cO_{\mP^2}(d)|$
we have a natural morphism: $\pi(\tSi)\to Gr(\mP^n,|\cO_{\mP^2}(d)|)$.
Hence there is a rational map: $\Aux-->Gr(\mP^n,|\cO_{\mP^2}(d)|)$, whose indeterminacy locus lies
in $\Aux\smin\pi(\tSi)$.
\\\parbox{9.2cm}{Resolve this map by a chain of blowups $\tAux\stackrel{\phi}{\to}\Aux$ (note that $\Aux$ is itself smooth),
then get the diagram on the right.

Now, the construction ensures that the projection $\ltSi'\to\tAux$
is a locally trivial fibration, whose fibres are linear subspaces of $|\cO_{\mP^2}(d)|$.
}\hspace{0.6cm}
\begin{picture}(0,0)(0,-25)
\put(0,0){$\tAux\times|\cO_{\mP^2}(d)|\supset\ltSi'\to\ltSi\supset\tSi\sset \Aux\times|\cO_{\mP^2}(d)|$}
\put(80,-10){$\da\hspace{0.85cm}\searrow\hspace{0.6cm}\searrow$}
\put(75,-25){$\tAux\stackrel{\phi}{\to} \Aux\supset\pi(\tSi)$}
\put(80,-35){$\da\hspace{1.8cm} \swarrow$}
\put(75,-50){$Gr\Big(\mP^n,|\cO_{\mP^2}(d)|\Big)$}
\end{picture}
\\
\\
This follows from the pull-back of the universal family on the Grassmanian.
In more detail, let $pt\in\tAux\smin\pi(\tSi)$. Let $(D,pt)\sset\tAux$ be the germ of a smooth curve such
that $D\smin\{pt\}\sset\pi(\tSi)$. The fibration over $D\smin\{pt\}$ is locally trivial and
extends (e.g. by the topological closure in the Grassmanian) to a locally trivial fibration over $D$.

Hence the fibre of $\ltSi'\to\tAux$ over $pt$ contains a linear subspace $L_{pt,D}\approx\mP^n\sset|\cO_{\mP^2}(d)|$.
But, this linear subspace does not depend on the curve $D$, it is determined just by the point
of the Grassmanian, the image of $\{pt\}$. Hence, by taking all the possible curves at $pt\in\tAux$ we get:
 all the fibres are linear subspace of $|\cO_{\mP^2}(d)|$ of constant dimension. Hence
the statement.
\epr
Note that the theorem is completely general. In practice it is quite hard to realize the resolved
base space $\tAux$ as some nice variety with a simple presentation of the (co)homology ring, explicit generators, etc.

\subsubsection{Flat limits and the decomposition of $\ltSi_{\sx\sy}|_{x=y}$  into
irreducible components}\label{Sec Background Collisions}
The following question occurs frequently. Given the singularity types $\sx\sy$
 find the defining equations of the equisingular stratum $\ltSi_{\sx\sy}(x,y,\dots )$  near the diagonal  $x=y$.
 And then
describe the decomposition of $\ltSi_{\sx\sy}|_{x=y}$  into irreducible components and find their defining ideals and
multiplicities. We call this process: collision.

In general this collision/adjacency problem is very complicated. It is manageable in the case of linear singularity types,
because the equations of
$\ltSi_\sx(x,\dots )$, $\ltSi_\sy(y\dots )$ are known (and hence the equations of $\ltSi_{\sx\sy}$ away from the diagonal).
To obtain the needed set of equations one expands $y=x+\ep v$ where $v\in l=\overline{xy}$,
and considers the ideal $\bl I_\sx(x,\dots ),I_\sy(x+\ep v,\dots )\br$.
Now one should take the flat limit as $\ep\to0$, a well known operation
in commutative algebra (cf. for example \cite[chapter XV]{Eisenbud-book}).
\\
\\
From the computational point of view one does the following.
Start from the ideal $\bl I_\sx(x,\dots ),I_\sy(x+\ep v,\dots )\br$, where in $I_\sy(x+\ep v,\dots )$ all the equations are series in $\ep$.
If for a member of this ideal the expansion in $\ep$ has no "constant" term, of 0'th order in $\ep$, divide it by the maximal
possible power of $\ep$ and add this normalized version to the ideal.
After several such iterations the process stabilizes,
i.e. no non-trivial syzygies remain.
So, we have obtained the defining ideal of $\ltSi_{\sx\sy}$ near the diagonal.
Then the ideal of the restriction $\ltSi_{\sx\sy}|_{x=y}$ is obtained by
putting $\ep=0$ in all the equations, i.e. constructing $I\otimes\cO_{x=y}=I/(\ep)$.

This process is described in details in \cite[\S 3.1]{Ker07Collisions} where many examples are considered.
\subsubsection{Collision with an \omp~ is the basic and most important case.}
\label{Sec Background Collisions With an OMP} Let $\sy$ denote the
singularity type of an \omp.
Assume $\mult(\sx)=p+1\ge \mult(\sy)=q+1$ and the collision is generic, i.e. the curve $l=\overline{xy}$ is
not tangent to any of the non-free branches of $\sx$, see definition \ref{DefNonFreeBranches}.

We should translate the conditions at the point $y$ to conditions at $x$.
Outside the diagonal $x=y$ the stratum is defined by the set of conditions corresponding to $\ltSi_{\sx}$, and by
the condition $f|_y^{(q)}=0$. The later is the (symmetric) form of derivatives of order $q$, calculated at the
point $y$ (in homogeneous coordinates). In the neighborhood of $x$ expand $y=x+\ep v$, here $\ep$ is
small and $v$ is the direction along the line $l=\overline{xy}$.

To take the flat limit, expand $f|_y^{(q)}$ around $x$, we get
$0=f|_y^{(q)}=f|_x^{(q)}+\dots +\frac{\ep^{p-q}}{(p-q)!}f|_x^{(p)}(\underbrace{v\dots v}_{p-q})+\dots $.
First several terms in the expansion vanish, up to the multiplicity of $\sx$.
Normalize by the common factor of $\ep$ to get the power series:
\beq\label{EqInitialSystemOrdMultPoin}
\frac{1}{(p-q+1)!}f|_x^{(p+1)}(\underbrace{v\dots v}_{p+1-q})+\frac{\ep}{(p-q+2)!} f|_x^{(p+2)}(\underbrace{v\dots v}_{p+2-q})+
\frac{\ep^2}{(p-q+3)!} f|_x^{(p+3)}(\underbrace{v\dots v}_{p+3-q})+\dots
\eeq
There are $\bin{q+2}{2}$ series here.
To take the flat limit, we should find all the syzygies between these series and the equations for $\ltSi_{\sx}$.
First we find the "internal" syzygies among the series themselves.
\bel
$\bullet$ The standard basis, obtained by considering all the syzygies of the equation (\ref{EqInitialSystemOrdMultPoin}), is
(with numerical coefficients omitted):
\beq\label{EqTriangularSystemOrdMultPoint}
\begin{tabular}{@{}>{$}c<{$}@{}>{$}c<{$}@{}>{$}c<{$}@{}>{$}c<{$}@{}>{$}c<{$}@{}>{$}c<{$}@{}>{$}c<{$}@{}>{$}c<{$}@{}>{$}c<{$}@{}>{$}c<{$}}
\\
 f|_x^{(p+1)}(\underbrace{v\dots v}_{p+1-q}) &
+& \ep f|_x^{(p+2)}(\underbrace{v\dots v}_{p+2-q})&
+&\ep^2f|_x^{(p+3)}(\underbrace{v\dots v}_{p+3-q})&
+&\ep^3f|_x^{(p+4)}(\underbrace{v\dots v}_{p+4-q})&+\dots
\\
0 &+&  f|_x^{(p+2)}(\underbrace{v\dots v}_{p+3-q})&
+&\ep f|_x^{(p+3)}(\underbrace{v\dots v}_{p+4-q})&
+&\ep^2 f|_x^{(p+4)}(\underbrace{v\dots v}_{p+5-q})&+\dots
\\
0 &+& 0&+& f|_x^{(p+3)}(\underbrace{v\dots v}_{p+5-q})
&+& \ep f|_x^{(p+4)}(\underbrace{v\dots v}_{p+6-q})&+\dots
\\
\dots &\dots &\dots &\dots &
\\
0&+&0 &+&0&+&\dots .&+&   f|_x^{(p+q+1)}(\underbrace{v\dots v}_{p+q+1})&+\dots
\end{tabular}
\eeq
\li If the singularity $\sx$ is an \omp~ ($x^{p+1}_1+x^{p+1}_2$, with $p\ge q$) then the series as
above define the stratum $\ltSi_{\sx\sy}(x,y,l)\sset \Aux\times|\cO_{\mP^2}(d)|$ near the diagonal $\{x=y\}$.
 Here $\Aux=\{(x,y,l)|x\in l\ni y\}\sset\mP^2_x\times\mP^2_y\times\tP^2_l$.
 In particular, the restriction $\ltSi_{\sx\sy}(x,y,l)|_{x=y}$ is an irreducible variety, whose generic point
corresponds to the singularity type of $(x^{p-q}_1+x^{p-q}_2)(x^{q+1}_1+x^{2q+2}_2)$ with $l$-the tangent
line to $(x^{q+1}_1+x^{2q+2}_2)$.
\eel
Note that all the series are linear in function/its derivatives in accordance
with proposition \ref{Thm Desingularization fo Strata}.
In several simplest cases we have: $A_1+A_1\ra A_3$, $D_4+A_1\ra D_6$, $X_9+A_1\ra X_{1,2}$, $D_4+D_4\ra J_{10}$,
 $X_9+D_4\ra Z_{13}$.
\\
\bpr $\bullet$
The syzygies are obtained as a consequence of the Euler identity for homogeneous polynomial
\mbox{$\sum x_i\di_if=\deg(f)f$}.
By successive contraction of the series of equation (\ref{EqInitialSystemOrdMultPoin}) and omitting $f^{(j)}|_x$
for $j\le p$ with $x$, we get the series
\beq\tinyA
\begin{tabular}{@{}>{$}c<{$}@{}>{$}c<{$}@{}>{$}c<{$}@{}>{$}c<{$}@{}>{$}c<{$}@{}>{$}c<{$}@{}>{$}c<{$}@{}>{$}c<{$}}
\\
 \frac{1}{(p-q+1)!}f|_x^{(p+1)}(\underbrace{v\dots v}_{p+1-q})&
+& \frac{\ep}{(p-q+2)!} f|_x^{(p+2)}(\underbrace{v\dots v}_{p+2-q})&
+& \frac{\ep^2}{(p-q+3)!} f|_x^{(p+3)}(\underbrace{v\dots v}_{p+3-q})+\dots &
\\
\frac{(d-p-2)}{(p-q+2)!}f|_x^{(p+1)}(\underbrace{v\dots v}_{p+2-q}) &
+& \frac{\ep (d-p-3)}{(p-q+3)!}f|_x^{(p+2)}(\underbrace{v\dots v}_{p+3-q})&
+& \frac{\ep^2(d-p-4)}{(p-q+4)!} f|_x^{(p+3)}(\underbrace{v\dots v}_{p+4-q})+\dots &
\\
\dots &\dots &\dots &\dots &
\\
\frac{\prod^{q+1}_{i=2}(d-p-i)}{(p+1)!}f|_x^{(p+1)}(\underbrace{v\dots v}_{p+1})&
+&\frac{\ep\prod^{q+1}_{i=2}(d-p-1-i)}{(p+2)!} f|_x^{(p+2)}(\underbrace{v\dots v}_{p+2})&
+&\frac{\ep^2\prod^{q+1}_{i=2}(d-p-2-i)}{(p+3)!} f|_x^{(p+3)}(\underbrace{v\dots v}_{p+3})+\dots &
\end{tabular}
\eeq
\\
Here the first row is the initial series, the second is obtained by contraction with $x$ once, the $p+2$'th
row is obtained by contracting ($p+1$) times with $x$.

Apply now the Gaussian elimination, to bring this system to the upper triangular form.
\ls Eliminate from the first column all the entries of the rows $2\dots (p+2)$. For this contract the first
row necessary number of times with $v$ (fix the numerical coefficient) and subtract.
\ls Eliminate from the second column all the entries of the rows $3\dots (p+2)$.
\ls \dots .
\\Normalize the rows (i.e. divide by the necessary power of $\ep$).

In this way we get the "upper triangular" system of series in eq. (\ref{EqTriangularSystemOrdMultPoint})
(we omit the numerical coefficients).
It remains to check that there are no more "internal" syzygies. This is proved in the second part.
\\
\\
\li Suppose $\sx$ is an \omp. First we should check that the upper triangular system, obtained above,
has no more syzygies. Alternatively, that it generates a prime ideal. This can be done directly
(by methods of commutative algebra). However in our case a simpler way is to check the dimension
of the restriction $\ltSi_{\sx\sy}|_{x=y}$ and to compare it to the dimension obtained from
(\ref{EqTriangularSystemOrdMultPoint}).

Note that $\dim(\Aux)=4$ and $\ltSi_{\sx\sy}$ is irreducible.
Hence the fibre over the generic point $(x,y,l)$ has dimension $\dim(\ltSi_{\sx\sy})-4$.
For the diagonal: $\dim\{(x,y,l)|x=y\in l\}=3$. Therefore
\beq
\dim(\ltSi_{\sx\sy}|_{x=y})\ge (\dim(\ltSi_{\sx\sy})-4)+3=\dim(\ltSi_{\sx\sy})-1
\eeq
On the other hand, by the irreducibility of $\ltSi_{\sx\sy}$ one has:
$\dim(\ltSi_{\sx\sy}|_{x=y})\le \dim(\ltSi_{\sx\sy})-1$. Thus we get the dimension of $\ltSi_{\sx\sy}|_{x=y}$
and of its fibre over the generic point $\{x=y\in l\}$.

Now check the scheme defined by equation (\ref{EqTriangularSystemOrdMultPoint}). Take the limit $\ep\to0$
(i.e. omit the higher order terms in each row):
\beq{
f|_x^{(p)}=0,~~~f|_x^{(p+1)}(\underbrace{v\dots v}_{p+1-q})=0,~~~f|_x^{(p+2)}(\underbrace{v\dots v}_{p+3-q})=0,
~~f|_x^{(p+3)}(\underbrace{v\dots v}_{p+5-q})=0~~\dots ,f|_x^{(p+q+1)}(\underbrace{v\dots v}_{p+q+1})=0
}
\eeq
\parbox{13.5cm}{
Observe that the system is linear in $f$. For a fixed $x,v$, e.g. $x=(0,0,1)$, $v=(0,1,0)$,
each equation is vanishing of some directional derivatives, and they are linearly independent.
Hence the fibres over $(x,v)$ are linear spaces of constant dimension and the variety is irreducible.

As each condition means the absence of some monomials in the expansion of $f$, the \ND~is easily constructed.
Now the co-dimension of the fibre over $\{x=y\in l\}$ can be computed e.g. as the number of $\mZ^2_{\ge0}$
points strictly below the \ND.
}
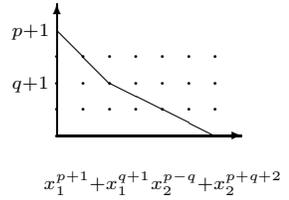
\begin{picture}(0,0)(-20,10)
\mesh{0}{0}{6}{3}{10}{70}{50}
\put(0,40){\line(1,-1){20}} \put(20,20){\line(2,-1){40}}
\put(-17,38){$\tinyA p+1$}  \put(-17,18){$\tinyA q+1$} 
\put(-5,-20){$\tinyA x^{p+1}_1+x^{q+1}_1x^{p-q}_2+x^{p+q+2}_2$}
\end{picture}
\\
For the dimension one gets precisely $\dim(\ltSi_{\sx\sy})-4$, i.e. the expected dimension.
This proves that the system in (\ref{EqTriangularSystemOrdMultPoint}) defines an irreducible variety:
the stratum $\ltSi_{\sx\sy}$ near the diagonal.
\\
\\
Finally, from the last set of equations or from the \ND~ we get: the restriction $\ltSi_{\sx\sy}|_{x=y}$
is irreducible. It is the equisingular stratum $\ltSi_\mS$ for the singularity type obtained
from the \ND.
Since all the slopes of the \D lie in the segment $[\frac{1}{2},2]$ the singularity type is linear.

For some fixed $x,l$ the generic point of $\ltSi_{\sx\sy}|_{x=y}$ corresponds to the generic function with the given
\ND, i.e. (with some numerical coefficients) $f=x^{p+1}_1+\cdots+x^{q+1}_1x^{p-q}_2+\cdots+x^{p+q+2}_2+$
{\it higher order terms}. By genericity $f$ is non-degenerate with respect to its diagram, hence the
singularity is the union of $(p-q)$ non-tangent smooth branches  and
$(q+1)$ smooth branches with pairwise tangency 2, i.e. has the singularity type of
$(x^{p-q}_1+x^{p-q}_2)(x^{q+1}_1+x^{2q+2}_2)$.
\epr
In the more general case, when $\sx$ is not an \omp, the generators of $I(\ltSi_{\sx})$ should be added
and one checks again for the possible syzygies.
\subsubsection{Tangency of the degenerating hypersurface and the diagonal.}\label{SecMultiplicityPiecesOverDiagonal}
Each step of the method is the intersection $\ltSi_1\cap V=\ltSi_2\cup R\subset|\cO_{\mP^2}(d)|\times \Aux$.
Here the ideal $I(\ltSi_1)$ is assumed to be known, $V$ is a hypersurface (with the known equation), $R$
is usually non-reduced, the multiplicities of its components are the degrees of tangency of $V$ and $\ltSi_1$
near the diagonal.

The multiplicities of $R$ are computed in a standard way.
The primary decomposition of the ideal $\bl I(\ltSi_1),I(V)\br$ consists of $I(\ltSi_2)$ and some
ideals of schemes over the diagonal $x=y$. These ideals give the multiplicities.

As an example, consider the degeneration of $\ltSi_{\sx\sy}$ by increasing the local degree of
intersection of the curve and the line $l=\overline{x}$ at $x\in\mP^2_x$. For the generic curve in
the original stratum this degree is just the multiplicity of $\sx$: $(p+1)$. Hence the degenerating hypersurface is
$V=\{(f|_x^{(p+1)})(\underbrace{v,\dots ,v}_{p+1})=0\}$.
 We assume that $v$  does not belong to the lines of the tangent cone at $x$, i.e. $l\neq l_x$.
Let $\ltSi_{\sx'\sy}$ be the degenerated stratum.
\bel\label{Thm Multiplicity Residual Piece Over Diagonal}
Suppose $\sx,\sy$ are \omp s of multiplicities $p+1\ge q+1$.
The multiplicity of the residual piece is $q+1$, i.e. $[\ltSi_{\sx\sy}]\times [V]=[\ltSi_{\sx'\sy}]+(q+1)[R_{reduced}]$
\eel
\bpr
The locally defining (prime) ideal of $\ltSi_{\sx\sy}$ near the diagonal $x=y$ is obtained in
\S\ref{Sec Background Collisions With an OMP}:
\beq
f|_x^{(p)}=0,~~~f|_x^{(p+1)}(\underbrace{v\dots v}_{p+1-q})+\dots =0,~~~f|_x^{(p+2)}(\underbrace{v\dots v}_{p+3-q})+\dots =0,
~~f|_x^{(p+3)}(\underbrace{v\dots v}_{p+5-q})+\dots =0~~\dots ,f|_x^{(p+1+q)}(\underbrace{v\dots v}_{p+1+q})+\dots =0
\eeq
Add the hypersurface equation $\{f|_x^{(p+1)}(\underbrace{v\dots v}_{p+1})=0\}$
 to this system of generators, decompose the so obtained ideal into primary components
 and find the resulting multiplicity.
 So, contract the first series $q$ times with $v$, the second series $(q-1)$ times with $v$ etc.:
\beq\ber
f|_x^{(p+1)}(\underbrace{v\dots v}_{p+1})+\ep^{q+1}f|_x^{(p+q+2)}(\underbrace{v\dots v}_{p+q+2})\dots =0,
~~~f|_x^{(p+2)}(\underbrace{v\dots v}_{p+2})+\ep^{p+q+2}f|_x^{(p+q+2)}(\underbrace{v\dots v}_{p+q+2})\dots =0,\\
~~f|_x^{(p+3)}(\underbrace{v\dots v}_{p+3})+\ep^{p+q+2}f|_x^{(p+q+2)}(\underbrace{v\dots v}_{p+q+2})\dots =0~~\dots ,
f|_x^{(p+q+1)}(\underbrace{v\dots v}_{p+q+1})+\dots =0
\eer\eeq
Together with $f|_x^{(p+1)}(\underbrace{v\dots v}_{p+1})=0$ the first equation becomes:
$\ep^{q+1}f|_x^{(p+q+2)}(\underbrace{v\dots v}_{p+q+2})\dots =0$.

So the ideal has two components:
$\bl \ep^{q+1},f|_x^{(p+1)}(\underbrace{v\dots v}_{p+1}),f|_x^{(p+2)}(\underbrace{v\dots v}_{p+2}),\dots ,f|_x^{(p+q+1)}(\underbrace{v\dots v}_{p+q+1})\br$
and
\beq\ber
\bl f|_x^{(p+1)}(\underbrace{v\dots v}_{p+1}),
~~~f|_x^{(p+2)}(\underbrace{v\dots v}_{p+2})+\ep^{p+q+2}f|_x^{(p+q+2)}(\underbrace{v\dots v}_{p+q+2})\dots ,\\
~~f|_x^{(p+3)}(\underbrace{v\dots v}_{p+3})+\ep^{p+q+2}f|_x^{(p+q+2)}(\underbrace{v\dots v}_{p+q+2})\dots ~~\dots ,
f|_x^{(p+q+1)}(\underbrace{v\dots v}_{p+q+1})+\dots ,f|_x^{(p+q+2)}(\underbrace{v\dots v}_{p+q+2})+\dots
\br
\eer\eeq
(In addition to these, the ideal contains also the equations of $\ltSi_{\sx\sy}$ near the diagonal.)

The first component is supported on the diagonal and is of multiplicity $q+1$, corresponding to the residual piece.
The second corresponds to $\ltSi_{\sx'\sy}$.
\epr
\section{The method}\label{Sec Method}
First we recall the case of curves with one singular point and state the necessary results.
Then we consider the case of two singular points and discuss additional complications.

Recall some conventions.
\li We work with many liftings, i.e. embeddings of $\Si_{**}$ into some multi-projective spaces.
To avoid messy notations
(like $\widetilde{\tSi}$) we denote by $\tSi$ {\it any} lifting. This causes no confusion, as we always specify the embedding
(e.g. $\tSi\subset|\cO_{\mP^2}(d)|\times \Aux$) or the associated parameters, e.g. $\tSi(x,y,l,\dots )$.
\li We always work with integral cohomology $H^*(\mP^n)=H^*(\mP^n,\mZ)=\mZ[t]/t^{n+1}$
or with integral homology.
\li
As was said above, once the cohomology class $[\ltSi_\mS]\in H^*(|\cO_{\mP^2}(d)|\times \Aux,\mZ)$ is computed,
the needed degree is obtained
by Gysin homomorphism, corresponding to the projection $\Aux\times|\cO_{\mP^2}(d)|\to|\cO_{\mP^2}(d)|$
(i.e. extraction of a particular coefficient).
If the projection $\tSi_\mS\to\Si_\mS$ is a covering, one should also divide by the order of the
symmetry group of branches.
Therefore, in the following we are interested in the cohomology classes of the {\it lifted}
strata $[\ltSi_{\sx\sy}]$.
\li Below we consider mostly linear singularities, for the definition and properties
cf.\S\ref{Sec Background Linear Singularities}.
\subsection{Equisingular strata of curves with one singular point}\label{Sec Method Strata of Unisingular Curves}
Here we summarize briefly some results of \cite{Ker06 enum.unising.curves}.
The equisingular strata are resolved by lifting to a bigger ambient space.
The lifted strata can be defined by some simple explicit conditions.
\subsubsection{The case of an \omp}\label{ExOneOrdinaryMultiplePoint}
The ordinary multiple point ($\mS=x^{p+1}_1+x^{p+1}_2$) is the simplest case (cf. \cite[\S2.3]{Ker06 enum.unising.curves}).
The first lifting (tracing the singular point) is already a {\it smooth, globally complete} intersection:
\beq
\ltSi_{\mS}(x)=\{(x,f)|~f|_x^{(p)}=0\}\subset\mP^2_x\times|\cO_{\mP^2}(d)|
\eeq
Here $f|_x^{(p)}$ is the tensor of all the partial derivatives of order $p$ in homogeneous coordinates,
calculated at the point $x$. This precisely encodes vanishing of the function and of all the derivatives up
to order $p$ in local coordinates, i.e. the conditions of an \omp.
The transversality of the defining conditions can be easily seen. For example, fix some point $x$,
then the conditions are linear, so the transversality is equivalent to linear independence.

Thus in this case the (co)homology class $[\ltSi_{\mS}(x)]$ is just the product of the (co)homology classes of
the defining hypersurfaces:
\beq
[\ltSi_{\mS}(x)]=\Big(F+(d-p)X\Big)^{\bin{p+2}{2}}\in H^*(\mP^2_x\times|\cO_{\mP^2}(d)|,\mZ)
\eeq
The class of $[\lSi_\mS]$ is obtained by Gysin projection, which in this case means:
to extract the coefficient of $X^2$. Hence
$[\lSi_\mS]=\bin{\bin{p+2}{2}}{2}(d-p)^2F^{\bin{p+2}{2}-2}\in H^*(|\cO_{\mP^2}(d)|,\mZ)$.
\subsubsection{The general linear singularity}
In general, let $\mS$ be a {\it linear} singularity type with the tangent cone $T_{\mS}=l_1^{p_1}\dots l_k^{p_k}$.
Here $p=\sum p_i$ is the multiplicity.
Consider the corresponding incidence variety, i.e. lift$\lSi_\mS$ to a bigger ambient space:
\beq
\ltSi_\mS(x,l_1\dots l_k):=
\overline{\{(x,l_1\dots l_k,f)|~C=f^{-1}(0)\sset\mP^2_x\text{ has the singularity of type $\mS$ at $x$
with $T_{(C,x)}=l_1^{p_1}\dots l_k^{p_k}$}\}}
\eeq
\parbox{13cm}
{The embeddings and projections of the varieties are given on the diagram. Here
\beq
\Aux=\{(x,l_1\dots l_k)| ~ \forall i ~ x\in l_i\}\sset\mP^2_x\times\tP^2_{l_1}\times\dots \times\tP^2_{l_k}
\eeq is the auxiliary incidence variety of lines through the points.
}~~~
$\bM \tSi_\mS\subset\ltSi_\mS\subset&|\cO_{\mP^2}(d)|\times \Aux\\\hspace{-0.2cm}\da\hspace{0.6cm}\da&\da\\
\Si_\mS\subset\lSi_\mS\subset&|\cO_{\mP^2}(d)|
\eM$
\\\\
Assume that the parts of the tangential decomposition of $\mS$ (cf.\S\ref{Sec Background SingularityTypes})
have distinct singularity types and are ordered. Then for a given $(C,0)$ the point $(x,l_1\dots l_k,C)\in\ltSi_\mS$
is unique.
Then the projection $\tSi_\mS\to\Si_\mS$ is an isomorphism and hence $\ltSi_\mS\to\lSi_\mS$ is birational.
If the singularity types for some of the parts coincide (e.g. for an \omp) then $\tSi_\mS\to\Si_\mS$
is an un-ramified covering.
\\
\\
As the type $\mS$ is linear, the fibres of the projection $\ltSi_\mS(x,l_1\dots l_k)\to \Aux$ are linear
 subspaces of $|\cO_{\mP^2}(d)|$. The cohomology class $[\ltSi_\mS]\in H^*(|\cO_{\mP^2}(d)|\times \Aux,\mZ)$
 can be computed as follows.

Let $\{f=0\}=(C,x)\sset(\mC^2,0)$ be the (generic) curve-germ of the type $\mS$.
Several initial jets of $f$ are divisible by the one-form $l_i$:
\beq
\jet_{p-1}(f)|_x=0, ~ ~ \jet_p(f)|_x=l^{p_i}_i\times (\dots .), ~ ~ \jet_{p+1}(f)|_x=l^{p_i-a_i}_i\times (\dots .),\dots .
 ~ ~ \jet_{p+n_i}(f)|_x=l_i\times (\dots .),
\eeq
To obtain this, choose the tangent line $l_i$
as one of the coordinate axes, the other axis is chosen generically.

These conditions can be written globally in $\Aux$ as:
\beq
f^{(p-1)}|_x=0, ~ ~f^{(p)}|_x\sim SYM(l^{p_i}_i,A_{p-p_i}), ~\dots . ~f^{(p+n_i)}|_x\sim SYM(l^{q_i}_i,A_{p+n_i-1})
\eeq
Here $A_j$ are some auxiliary symmetric tensors (of rank $j$) and $SYM$ means the symmetrization.
Do the same for all the tangent lines to get the collection of conditions:
\beq
f^{(p-1)}|_x=0, ~ ~f^{(p)}|_x\sim SYM(l^{p_1}_1,\dots .,l^{p_k}_k), ~\dots .
~f^{(p+n)}|_x\sim SYM(l^{n_1}_1,\dots ,l^{n_k}_k,A_{p+n-n_1-\dots -n_k})
\eeq
Now use the Euler identity $f^{(N)}|_x(x)\sim f^{(N-1)}|_x$ and its consequences to present these conditions
in the form:
\beq\label{Eq Defining Eqs Equis Stratum One Sing Point}
f^{(p+n)}|_x\sim SYM(l^{n_1}_1,\dots ,l^{n_k}_k,A_{p+n-n_1-\dots -n_k}), ~ ~
A_{p+n-n_1-\dots -n_k}(x)\sim SYM(l^{m_1}_1,\dots ,l^{m_k}_k,A_{p+n-\sum n_i-\sum m_i}), ~ \dots ., ~
A_{q}(x)=0
\eeq
Now observe that distinct proportionality conditions are obviously mutually transversal. For example
 $f$ appears in the first proportionality only, $A_{p+n-n_1-\dots -n_k}$ only in the first and in the second etc.
\bthe\cite[\S 3]{Ker06 enum.unising.curves}
Let $\ltSi_\mS(x,\{l_i\},\{A_j\})$ be the lifted stratum, defined by the
equation (\ref{Eq Defining Eqs Equis Stratum One Sing Point}).
Then $\ltSi_\mS(x,\{l_i\},\{A_j\})\to \Aux$ is the (smooth) projective bundle over the auxiliary incidence variety
 $\Aux$ of points, lines $\{l_i\}$ and symmetric tensors $\{A_{**}\}$. The fibres are linear subspaces of
 $|\cO_{\mP^2}(d)|$. The cohomology class of the lifted stratum is the product of the cohomology
 classes
\beq\ber
[\ltSi_\mS(x,\{l_i\},\{A_j\})]=\Big[f^{(p+n)}|_x\sim SYM(l^{n_1}_1,\dots ,l^{n_k}_k,A_{p+n-n_1-\dots -n_k})\Big]\times
\\
\Big[A_{p+n-n_1-\dots -n_k}(x)\sim SYM(l^{m_1}_1,\dots ,l^{m_k}_k,A_{p+n-\sum n_i-\sum m_i})\Big]
 \times\dots \times\Big[A_{q}(x)=0\Big]
\eer\eeq
The projection  $\tSi_\mS(x,\{l_i\},\{A_j\})\to\tSi_\mS(x,\{l_i\})\to\Si_\mS$ is isomorphism
(or unramified cover if some tangential components are equisingular).
\ethe
Note that each of the conditions is proportionality of two tensors, i.e. coincidence of two points in some
big projective space. Correspondingly its cohomology class is readily written,
from \S\ref{Sec Background Blowup Diagonal}.
\bex\label{Sec Method Examples One Sing Point Many Branches}
Let $\mS$ be the type of $k$ pairwise non-tangent branches, each of the cuspidal type $x^{p_\al}_1+x^{p_\al+1}_2$.
For example, a representative of this type can be chosen as: $f=\prod l^{p_\al}_\al+f_{p+1}$ where $f_{p+1}$
is generic, and $\{l_\al\}$ are some distinct linear forms.
Then $T_{\mS}=l^{p_1}_1\dots l^{p_k}_k$ and any other germ $g=g_p+g_{p+1}+\dots $ with the same tangent cone and
generic $g_{p+1}$ is of this type. Hence the lifted stratum is
\beq\ber
\ltSi_\mS(x,l_1\dots l_k)=\Big\{(x,l_1\dots l_k,f)|~ f^{(p)}\sim SYM(l^{p_1}_1\dots l^{p_k}_k),~~\forall \al~x\in l_\al\Big\}
\sset \Aux\times|\cO_{\mP^2}(d)|
\\
\Aux=\{(x,l_1\dots l_k)|~\forall \al:~x\in l_\al\}\subset\mP^2_x\times\tP^2_{l_1}\times\dots \times\tP^2_{l_k}
\eer\eeq
Thus the cohomology class of the lifted stratum is:
\beq
[\ltSi_\mS(x,l_1\dots l_k)]=\prod^k_{\al=0}(X+L_\al)\sum^{\bin{p+2}{2}-1}_{j=0}(F+(d-p)X)^{\bin{p+2}{2}-1-j}(\sum p_\al L_\al)^j
\eeq
The class of the original stratum $[\lSi_\mS]$ is obtained by extracting the coefficient of $X^2L^2_1\dots L^2_k$.
The numerical formula is given in Appendix.
If some of $p_\al$ coincide one should also divide by the corresponding symmetry factor.
\eex
\subsubsection{Some degenerating divisors in the stratum $\ltSi_\mS(x,l)$}
\label{Sec Method CohomologyClassDegenerations}~\\
\parbox{14.5cm}
{We need the class of a specific divisor in $\ltSi_\mS$. Suppose the generic curve of the type $\mS$ can
be brought, in local coordinates, to some fixed \ND~$\Ga$. Let $(p,q)\in\Ga$
be a particular vertex, i.e. a point of $\mZ^2_{\ge0}\cap\Ga$, not in the interior of an edge of $\Ga$.
Erase this vertex, let $\Ga'=Conv(\Ga\smin(p,q))$ be the obtained diagram. Consider the
substratum $\Si_2\sset\Si_1=\Si_\mS$ parameterizing of all the curves that can be brought to the diagram $\Ga'$.
Assume $\Si_2\subsetneqq\Si_1$ so that this is a genuine degeneration. As $d$ is high
$\lSi_2$ is a hypersurface in $\lSi_1$. Its cohomology class can be computed as
follows.
}
\begin{picture}(0,0)(-10,20)
\put(0,0){\vector(1,0){60}}\put(0,0){\vector(0,1){60}}
\put(0,45){\line(1,-1){25}}  \put(25,20){\line(3,-1){30}} \put(22,17){$\bullet$}\put(12,28){$\bullet$}
\put(51,8){$\bullet$}\put(52,11){\line(4,-1){30}}
\multiput(14,31)(16,-8){3}{\line(2,-1){10}}
\put(-8,17){$p$}  \put(22,-8){$q$}
\put(4,55){$\hat{x}_2$}  \put(65,-5){$\hat{x}_1$}
\end{picture}
\bel\cite[section A.1.2]{Ker06 enum.unising.curves}\label{Thm Class of Degenerating Divisor}
The classes satisfy:
\beq
[\ltSi_2(x,l)]\Big(F+(d-p-2q)X+(q-p)L\Big)=[\ltSi_1(x,l)]\in H^*(\mXL\times|\cO_{\mP^2}(d)|)
\eeq
\eel
In particular, the class of $\ltSi_2(x,l)$ in $\ltSi_1(x,l)$ is the pull-back:
$i^*(F)+(d-p-2q)i^*(X)+(q-p)i^*(L)$.
\\
\bpr Let $l$ be the tangent line corresponding to the $\hat{x}_1$ axis.
We need to consider an additional point $v\in l$ with $v\neq x$. Correspondingly lift the strata to a bigger
auxiliary space $\mP^2_x\times\mP^2_v\times\tP^2_l\times|\cO_{\mP^2}(d)|$ as follows:
\beq
\ltSi_j(x,v,l):=\{(x,v,l,f)|~~(x,l,f)\in\ltSi_j(x,l),~~v\in l\}
\eeq
For the cohomology classes one has: $[\ltSi_j(x,v,l)]=[\ltSi_j(x,l)]\times[v\in l]$.

Observe that the line $l$ is spanned by $x,v$, provided $x\neq v$. Hence the degenerating
condition, erasing the vertex $(p,q)$, is: the directional derivative vanishes
$f^{(p+q)}|_x(\underbrace{v\dots v}_{q})_{\underbrace{0\dots 0}_{p}}=0$.
Here we take a particular component of the tensor. Alternatively we could consider
$f^{(p+q)}|_x(\underbrace{v\dots v}_{q}{\underbrace{\tv\dots \tv}_{p}})=0$ for some fixed generic vector $\tv$.

This condition is transversal to the stratum $\ltSi_j(x,v,l)$ provided: $v\neq x$ and the point $(1,0,0)$
does not lie on the line $l$, i.e. $l^p_0\neq0$.

Thus for the cohomology classes:
\beq
[\ltSi_2(x,l)]\times\Bigg([v\in l]\Big([f^{(p+q)}|_x(\underbrace{v\dots v}_{q})_{\underbrace{0\dots 0}_{p}}=0]-[l^p_0=0]\Big)-
multiplicity\times[x=v]\Bigg)=[\ltSi_1(x,l)]\times[v\in l]
\eeq
Here the diagonal $[x=v]$ is subtracted with some multiplicity, corresponding to the tangency
of the degeneration. The multiplicity can be computed directly, but it is easier to obtain as follows.
Note that the point $v$ moves freely along the line $l$, i.e. all the participating varieties
are fibrations with the fibre $\mP^1_v$. Hence in the cohomology classes all the terms quadratic in $V$ must cancel.
This fixes the multiplicity uniquely: $mult=q$.

 Finally, substitute to the formula the expressions for the classes
\beq
[v\in l]=V+L, ~ ~ [x=v]=X^2+XV+V^2, ~ ~ ~[f^{(p+q)}|_x(\underbrace{v\dots v}_{q})_{\underbrace{0\dots 0}_{p}}=0]=F+(d-p-q)X+qV
\eeq
From the obtained equation extract the part proportional to $V$, this corresponds to the projection
$\mP^2_x\times\mP^2_v\times\tP^2_l\times|\cO_{\mP^2}(d)|\to\mP^2_x\times\tP^2_l\times|\cO_{\mP^2}(d)|$.
Then one gets the equation of the lemma.
\epr
\subsubsection{Killing tangent cone.}\label{Sec Method Killing Tangent Cone}
A particular kind of degeneration that happens to be especially useful is: to increase the multiplicity of $\mS$ by one.
Given a singularity type $\mS$ of multiplicity $p$, let $\lSi_1\sset\lSi_\mS$ be
the subvariety of curves with multiplicity at least $p+1$.
In general $\lSi_1$ is reducible and does not correspond to a specific singularity type.
For example $\lSi_{E_6}\cup \lSi_{D_6}\sset\lSi_{A_5}$.
If $\mS$ is a linear type, $\lSi_1$ is much more restricted.
\bdp
If $\mS$ is linear then the subvariety as above, $\lSi_1\sset\lSi_\mS$, is an irreducible equisingular stratum of a
unique linear type $\mS'$. The degeneration $\mS\ra\mS'$ is called:
{\bf killing the tangent cone}.
\edp
\bpr
Consider the defining ideal of the lifted stratum $\ltSi_\mS$ as in equation
(\ref{Eq Defining Eqs Equis Stratum One Sing Point}).
By definition the ideal of $\lSi_1$ is obtained
by adding to the equations of $\ltSi_\mS$ the equations $f^{(p)}|_x=0$.

This bigger set of equations is still linear in the coefficients of $f$ and the ideal is obviously prime.
Hence the generic point of $\lSi_1$ is well defined its singularity type too.
From the defining set of equations we get that $\mS'$ is a linear type.
\epr
The cohomology class of the degeneration is given by:
\bel\label{Thm Kill Tangent Cone Cohom Class} Let $T_\mS=l_1^{p_1}\dots k_k^{p_k}$ be the tangent cone.
For the lifting $\ltSi_\mS(x,\{l_i\})$, the subvariety $\ltSi_{\mS'}\sset\ltSi_\mS$ is a divisor and its cohomology class is:
\beq
[\ltSi_{\mS}(x,\{l_i\})]\Big(F+(d-p)X-\sum p_iL_i\Big)=[\ltSi_{\mS'}(x,\{l_i\})]\in H^*(\Aux\times|\cO_{\mP^2}(d)|)
\eeq
\eel
In the paper we denote this divisor by $\kill T_\mS$.
\\
\bpr
Intersect the lifted stratum
with the hypersurface $\{(f|_x^{(p)})_{\underbrace{0\dots .0}_{p}}=0\}$. This means that the particular entry of
the tensor of derivatives vanishes.
(Alternatively, one could take any collection of fixed vectors $\{v_i\}$ $\{(f|_x^{(p)})(v_1\dots .v_p)=0\}$.)

Outside the locus $(l^{p_1}_1\dots l^{p_r}_r)_{0\dots 0}=0$ this kills the tangent cone. Therefore we get the equation for the cohomology classes:
\beq
[\ltSi_{\mS}(x,\{l_i\})]\Big([(f|_x^{(p)})_{0\dots 0}=0]-[(l^{p_1}_1\dots l^{p_k}_k)_{0\dots 0}=0]\Big)=[\ltSi_{\mS'}(x,\{l_i\})]
\eeq
which proves the statement.
\epr

\subsubsection{Rectifying the branch.}This degeneration is useful in the degeneration of the curve to a reducible
one (chipping off a line).
\\
\parbox{14cm}
{Consider a singular germ $(C,0)$ with $l$ one of the lines of the tangent cone. The goal is to increase
the degree of the intersection: $\deg(C\cap l,0)=k$ to $(k+1)$.
To do this, choose $l$ as one of the coordinate axes, cf. the \ND.
}
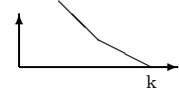
\begin{picture}(0,0)(-10,0)
\put(0,0){\vector(1,0){60}}\put(0,0){\vector(0,1){20}}\put(50,0){\line(-2,1){20}}\put(30,10){\line(-1,1){15}}
\put(48,-8){\tinyT k}
\end{picture}

So, the degeneration is: the monomial $x^k_2$ should be absent. The corresponding cohomology class is given
in \S\ref{Sec Method CohomologyClassDegenerations}, in our case it is: $F+(d-2k)X+kL$.
\subsubsection{Another way to compute the class $[\ltSi_\mS]$}\label{Sec Method One Sing Point Stepwise Intersections}
Using the class of the degenerating divisors we can compute the class $[\ltSi_\mS]$ in another way:
by a chain of intersections with hypersurfaces.

Let $p=\mult(\mS)$ and $T_\mS=l^{p_1}_1\dots l^{p_k}_k$.
Let $x\neq v_\al\in l_\al$ be some points on the lines. Then the defining
conditions (\ref{Eq Defining Eqs Equis Stratum One Sing Point}) can be formulated as vanishing of particular
entries of the derivative-tensors: $\Big\{ f^{(n)}(v_{i_1}\dots v_{i_n})=0\Big\}$, each means the absence
of a particular monomial corresponding to some point on the \ND. Then the class $[\ltSi_\mS]$
is obtained as the product of the classes of these divisors.
\bex Consider the cusp $\mS=x^p_1+x^{p+1}_2$.
Locally a curve-germ with such a singularity satisfies (we assume the $\hat{x}_2$ axis to be the tangent line):
\beq
\jet_{p-1}f=0, \di^p_2f=0, ~ \di^{p-1}_2\di_1f=0, ~ \dots .\di_2\di^{p-1}_1f=0
\eeq
So the stratum is obtained from the stratum of curves with an \omp~ (defined by $f^{(p-1)}|_x=0$)
 by intersection with $p$ degenerating hypersurface.
The cohomology class of each such hypersurface is given in proposition \ref{Thm Class of Degenerating Divisor}.
In total we have:
\beq
[\ltSi_\mS(x,l)]=(F+(d-p+1)X)^{\bin{p+1}{2}}\prod_{i=0}^{p-1}(F+(d+i-2p)X+(p-2i)L)
\eeq
\eex
\subsection{Equisingular strata of curves with two singular points}\label{Sec Method Strata For two singular points}
~\\
\parbox{11cm}
{Assume $\mult(\sx)\ge \mult(\sy)$. As in the case of one singular point, start from the lifting
$\tSi_{\sx\sy}\subset|\cO_{\mP^2}(d)|\times \Aux_\sx\times \Aux_\sy$.
To simplify the computation it is often useful to lift further ($\tSi_{\sx\sy}\subset|\cO_{\mP^2}(d)|\times \Aux$), taking into account various parameters
relating the two singular points, for example the line $l=\overline{xy}$.
If the types $\sx,\sy$ are distinct and the types of tangential components of $\sx$ are distinct (and the same for $\sy$) then
the projection $|\cO_{\mP^2}(d)|\times \Aux\to |\cO_{\mP^2}(d)|$ restricts to
the birational morphism $\ltSi_{\sx\sy}\to\lSi_{\sx\sy}$.
}~~~
$\bM \tSi_{\sx\sy}\subset\ltSi_{\sx\sy}\subset&|\cO_{\mP^2}(d)|\times \Aux\\\hspace{-0.4cm}\da\hspace{1cm}\da&\da\pi
\\\tSi_{\sx\sy}\subset\ltSi_{\sx\sy}\subset&|\cO_{\mP^2}(d)|\times \Aux_\sx\times \Aux_\sy
\\\hspace{-0.4cm}\da\hspace{1cm}\da&\da
\\\Si_{\sx\sy}\subset\lSi_{\sx\sy}\subset&|\cO_{\mP^2}(d)|
\eM$

\bex\label{Ex Lifting For 2OMP}
In the simplest case of  two \omp s there is only one (natural) additional parameter: the line $l=\overline{xy}$.
Note that for the generic configuration of the singularities (i.e. $x\ne y$) the line is already fixed, while
 for $x=y$ the line varies in a pencil.
 Therefore this additional
lifting is the embedded blow-up over the diagonal $\{x=y\}\subset\mP^2_x\times\mP^2_y$. Explicitly,
let $\sx=x^{p+1}_1+x^{p+1}_2,$ $\sy=y^{q+1}_1+y^{q+1}_2$ with $p>q$, then the lifting is defined by:
\beq\ber
\ltSi_{\sx\sy}(x,y,l)=\overline{\Big\{\ber(x,y,l,f)\\x\ne y\eer|~f|_x^{(p)}=0=f|_y^{(q)},~~~
l=\overline{xy}\Big\}}\subset\ \Aux\times|\cO_{\mP^2}(d)|,\\
\Aux:=\{(x,y,l)|~~~x\in l\ni y\}\subset\mP^2_x\times\mP^2_y\times\tP^2_l
\eer\eeq
Outside the diagonal $\{x=y\}$, the vanishing of the derivatives
is precisely the condition of \omp s. The fibres over the diagonal correspond to the singularity
obtained by "collision" of the two \omp s, the defining equations and the singularity type are obtained by the flat
limit of the equations $f|_x^{(p)}=0=f|_y^{(q)}$, cf.\S\ref{Sec Background Collisions With an OMP}.
\eex
\subsubsection{The direct approach: Stepwise intersection with hypersurfaces}\label{Sec Method Stepwise Intersection}
For linear type $\mS$ the defining ideal and the cohomology class of $\ltSi_\mS$ were given in
\S\ref{Sec Method Strata of Unisingular Curves}.
The stratum of curves with two singularities is the
closure: $\overline{\ltSi_\sx\underset{x\neq y}{\cap}\ltSi_\sy}$.
The naive intersection results in a non-pure dimensional scheme:
\beq
\ltSi_\sx{\cap}\ltSi_\sy=\ltSi_{\sx\sy}\cup R_{x=y}
\eeq
Here the "residual piece"(the contribution from the diagonal) $R_{x=y}$ is typically non-reduced and
of dimension higher than that of the needed stratum.
\\
\\
A way to repair this situation is to split the intersection into a step-by-step procedure of intersection with
hypersurfaces (cf. \cite{StuVog82}, \cite{vGas89}). 
Start from $\ltSi_\sx$ whose class is known.

Let $\ltSi_{\sx}\subset|\cO_{\mP^2}(d)|\times \Aux_x$ be the lifting as described above.
Corresponding to the lifting $\ltSi_{\sx\sy}\subset|\cO_{\mP^2}(d)|\times \Aux$ of
the diagram above and the projection $\pi$, define $\ltSi^{(y)}_\sx:=\pi^{-1}(\ltSi_{\sx}\times \Aux_y)\subset|\cO_{\mP^2}(d)|\times \Aux$.
For $\sx$ linear, the stratum $\ltSi_{\sx}$ is smooth and by the properties of $\pi$: $\ltSi^{(y)}_\sx$ is smooth too.

The cohomology class $[\ltSi^{(y)}_\sx]\in H^*(\Aux\times|\cO_{\mP^2}(d)|)$ is just the pull-back of $[\ltSi_\sx]$.
Consider the stratum $\ltSi_{\sx\sy}$ as a subvariety:
$\ltSi_{\sx\sy}\subset\ltSi^{(y)}_\sx\subset|\cO_{\mP^2}(d)|\times \Aux$.
Correspondingly, it is enough to calculate its class:
 $[\ltSi_{\sx\sy}]\in H^*(\ltSi^{(y)}_\sx,\mZ)$. The resulting class in  $H^*(|\cO_{\mP^2}(d)|\times \Aux,\mZ)$
 is then obtained by the pushforward.
So, $\ltSi^{(y)}_\sx$ is the starting point and we reach the stratum $\ltSi_{\sx\sy}$ by
successive intersections as follows.

As was explained in \S\ref{Sec Method One Sing Point Stepwise Intersections} the stratum $\ltSi_{\sy}$
can be obtained as a chain of degenerations:
$\ltSi_\sy=\ltSi_k\sset\ltSi_{k-1}\sset\dots \sset\ltSi_0=|\cO_{\mP^2}(d)|\times \Aux_y$.
Here at each step $\ltSi_j$ is a divisor in $\ltSi_{j-1}$, whose class is known.
Take the total preimages of this chain in $|\cO_{\mP^2}(d)|\times \Aux$, denote them by the same letters.
Intersect $\ltSi^{(y)}_\sx$ with these subvarieties one-by-one.
\beq\ber\label{Eq Stepwise Intersection With Hypersurfaces}
\ltSi^{(y)}_{\mS_x}\cap \ltSi_{1}=\ltSi_{\mS_x\mS_1}\cup R_{1},\ldots\ldots
\ltSi_{\mS_x\mS_j}\cap \ltSi_{j+1}=\ltSi_{\mS_x\mS_{j+1}}\cup R_{j+1},
\\
\ltSi_{\mS_x\mS_y}\subsetneq\ltSi_{\mS_x\mS_{k-1}}\subsetneq\ltSi_{\mS_x\mS_{k-2}}\subsetneq\dots
\subsetneq\ltSi_{\mS_x\mS_1}\subsetneq\ltSi_{\mS_x}^{(y)}
\eer\eeq
Here $\ltSi_{\mS_x\mS_{j+1}}$ is the needed piece (defined as the closure:
$\overline{\ltSi_{\mS_x\mS_j}\underset{x\neq y}{\cap} \ltSi_{j+1}}$).
Its dimension drops precisely by one with each intersection. $R_{j}$ is the
residual piece produced at the $j$'th step (it contains all
the non-enumerative contributions).

\bprop
$\bullet$
Each intersection $\ltSi_{\mS_x\mS_j}\cap \ltSi_{j+1}$ results in a pure
dimensional scheme  $\ltSi_{\mS_x\mS_{j+1}}\cup R_{j+1}$.
\li For linear types $\sx\sy$, the scheme $\overline{\ltSi^{(y)}_\sx\underset{x\neq y}{\cap^{j+1}_{i=1}}\ltSi_i}$
is irreducible and reduced.
\eprop
The residual piece $R_{j+1}$ is typically reducible and non-reduced.\\
\bpr $\bullet$
Note that the initial variety is irreducible, thus each components of the resulting
union $\overline{\ltSi^{(y)}_\sx\underset{x\neq y}{\cap^j_{i=1}}\ltSi_i}\cap \ltSi_{j+1}$ is of (strictly) smaller
dimension.
Conversely, the dimension of each component cannot drop by more than one
(since the intersection is with a hypersurface and the ambient space is irreducible).
\li Note that the defining equations of both $\overline{\ltSi^{(y)}_\sx\underset{x\neq y}{\cap^j_{i=1}}\ltSi_i}$
and $\ltSi_{j+1}$ are linear in $f$
(i.e. in the coordinates of $|\cO_{\mP^2}(d)|$).
\epr
Thus the contribution of the residual piece can be subtracted:
$[\ltSi_{\mS_x\cap^{j+1}_{i=1}\ltSi_i}]$$=[\ltSi_{\mS_x\cap^j_{i=1}\ltSi_i}][\ltSi_{j+1}]$$-[R_{j+1}]$ (with
multiplicities for $R_{j+1}$).
By repeating this procedure we calculate the needed class $[\ltSi_{\mS_x\mS_y}]$.
\\
\\
Therefore the enumerative problem is reduced to geometry of the residual pieces and their multiplicities.
This method is applied to the case of two \omp s in \S\ref{Sec Example Two OMP's}.
Explicit numerical results are given in the theorem
\ref{ThmTwoOrdinaryMultiplePoints}, corollary \ref{ThmTwoOrdinaryMultPointNumericalCorollary} and
in the Appendix \ref{AppendixNumericalResults}.
\subsubsection{Degenerations to higher singularities}\label{Sec Method Degenerations To OMP}
Let $\sx\sx'$ be two adjacent singularity types, i.e. $\lSi_{\sx'}\subsetneq\lSi_\sx\subset|\cO_{\mP^2}(d)|$.
Let $Degen\subset|\cO_{\mP^2}(d)|\times \Aux_x$ be a (reduced) "degenerating cycle" such that
$Degen\cap\ltSi_\sx=\ltSi_{\sx'}$ (at least set theoretically).
These degenerating cycles and their cohomology classes are described in \S\ref{Sec Method CohomologyClassDegenerations}.
Then for the cohomology classes we have (with some multiplicities): $[\ltSi_\sx][Degen]\sim[\ltSi_{\sx'}]\in H^*(|\cO_{\mP^2}(d)|\times \Aux)$.
The multiplicity here is the local degree of intersection of $Degen$ with $\ltSi_\sx$.

This equation allows to calculate $[\ltSi_{\sx'}]$ in terms of $[\ltSi_\sx]$.
The key observation is that such a degeneration is always "invertible":
the above equation has unique solution for $[\ltSi_\sx]$ in terms of $[\ltSi_{\sx'}]$ and $[Degen]$.
This happens because the class $[Degen]$ contains a monomial: a power of the generator of $H^*(|\cO_{\mP^2}(d)|,\mZ)$.
And $d$ is arbitrary high, so this monomial is not a zero divisor when multiplying bounded classes, i.e. lying
in $\oplus^k_{i=0}H^i(|\cO_{\mP^2}(d)|)$ for $k\ll\dim|\cO_{\mP^2}(d)|$  (cf. \cite[\S 2.2]{Ker06 enum.unising.curves}).
\\
\\
Using this, we can degenerate to some higher singularity types for which the classes are known.
In \S\ref{Sec Intro Example 2OMP Degen to Reducible} we degenerate a curve with two \omp s to a reducible curve,
with the line through the points as a component.
Once the cohomology class of $\ltSi_{x^{p+1}_1+x^{p+1}_2,x^{q+1}_1+x^{q+1}_2}$ is obtained,
we can continue to other singularities, degenerating them to a pair of \omp s.
\\
\\
A suitable degeneration is by increasing the multiplicity of $\sx$ or $\sy$. For the
original stratum $\lSi_{\sx\sy}\sset|\cO_{\mP^2}(d)|$ such degeneration is usually of high codimension. To turn this into
a degeneration in codimension 1, i.e. intersection by a hypersurface, consider the lifting that traces
all the lines of the tangent cone: $\ltSi_{\sx\sy}(x,y,l,\{l_{x,i}\},\{l_{y,j}\})$.
Then the degeneration is done by intersection with $\kill T_\mS$ divisor from \S\ref{Sec Method Killing Tangent Cone}.

Hence the algorithm consists of several steps: each time increasing the multiplicity of $\sx$ or $\sy$.
 At each step the variety is intersected by
the degenerating hypersurface, the result is reducible:
\beq
\ltSi_{\mS_x\mS_y}\cap \kill T_\sx=\ltSi_{\sx'\sy}\cup R^{x=y}
~ ~ ~ \text{ or } \ltSi_{\mS_x\mS_y}\cap \kill T_\sy=\ltSi_{\sx\sy'}\cup R^{x=y},
\eeq
\bcor
Let $\sx\sy$ be linear types. The degeneration approach reduces the enumerative problem to understanding
the irreducible components of $\ltSi_{\sx\sy}|_{x=y}$.
\ecor
\bpr
First note that each step gives an equation (in the cohomology ring of the ambient space $\Aux\times|\cO_{\mP^2}(d)|$)
from which the class of $\ltSi_{\mS_x\mS_y}$
 is restored in terms of the classes of
 $\ltSi_{\mS_x'\mS_y}$, $\kill T_\sx$ (or $\ltSi_{\mS_x\mS_y'}$, $\kill T_\sy$), $R^{x=y}$.

Indeed, the resulting variety $\ltSi_{\mS_x'\mS_y}\cup R^{x=y}$ is pure dimensional.
The class of $\kill T_\sx$ is given in \S\ref{Sec Method Killing Tangent Cone}.
The residual pieces $R^{x=y}$ are (as sets)
unions of the equisingular strata $\ltSi_{\mS}$, where ${\mS}$'s are the results of collision \ref{Sec Background Collisions}.
In other words: $\ltSi_{\sx\sy}|_{x=y}=\cup \ltSi_{\mS}$.
The way to compute the multiplicities of the decomposition $[\ltSi_{\sx\sy}|_{x=y}]=\sum m_i[\ltSi_{\mS_i}]$ is
 explained in \S\ref{SecMultiplicityPiecesOverDiagonal}.

Finally, after a finite number of steps the algorithm stops. This statement is immediate since e.g.
one can reach an \omp~ in a
finite number of steps. In fact, by killing the tangent cone, cf. \S\ref{Sec Method Killing Tangent Cone} one can
do this in less than $o.d.(\mS)-\mult(\mS)$ steps.
\epr
Examples of the degenerations are in \S\ref{Sec Example Degeneration to OMP's}. Below we discuss the
relevant components of $\ltSi_{\sx\sy}|_{x=y}$ contributing to residual pieces.
\subsubsection{Residual pieces for the degeneration $\ltSi_{\sx\sy}\to\ltSi_{\sx'\sy}$.}\label{Sec Residual Pieces sx to OMP}
Assume $\sx,\sy$ are linear types and $\mult(\sx)\ge \mult(\sy)$. Let $\kill T_\sx$ denote the "killing tangent cone" divisor
(as in lemma \ref{Thm Kill Tangent Cone Cohom Class}), so that $\ltSi_\sx\cap \kill T_\sx=\ltSi_{\sx'}$. For the case
of two singular points:
\beq
\ltSi_{\sx\sy}\cap \kill T_\sx=\ltSi_{\sx'\sy}\cup\Big(\ltSi_{\sx\sy}|_{x=y}\cap \kill T_\sx\Big)
\eeq
Let $\ltSi_{\sx\sy}|_{x=y}=\cup_j \ltSi_j$ be the decomposition into irreducible components.
As $\dim(\ltSi_j)\le \dim(\ltSi_{\sx\sy})-1$ the component $\ltSi_j$ contributes to the
equation in cohomologies iff $\ltSi_j\sset \kill T_\sx$.

From the local defining equation of $\kill T_\sx$ (\S\ref{Sec Method Killing Tangent Cone}) we get therefore:
\bprop
\beq
[\ltSi_{\sx\sy}][\kill T_\sx]=[\ltSi_{\sx'\sy}]+\sum_{j\in J}tang_j[\ltSi_j]
\eeq
where:
\ls $J$ indexes the irreducible components of $\ltSi_{\sx\sy}|_{x=y}$ corresponding to
collisions $\sx+\sy\to\mS$ with $\mult(\mS)>\mult(\sx)$.
\ls $\ltSi_j$ is the corresponding stratum lifted to the given ambient space.
\ls $tang_j$ is the degree of tangency of $\kill T_\sx$ and $\ltSi_{\sx\sy}$ along $\ltSi_j$.
\eprop
\subsubsection{Residual pieces for the degeneration $\ltSi_{\sx\sy}\to\ltSi_{\sx\sy'}$.}
Now the degeneration is done by intersection with the divisor $\kill T_\sy$. As $\mult(\sy)\le \mult(\sx)$
all the irreducible components of $\ltSi_{\sx\sy}|_{x=y}$ will contribute to the cohomology class:
\beq
[\ltSi_{\sx\sy}][\kill T_\sy]=[\ltSi_{\sx\sy'}]+\sum_{\ber components~of\\\ltSi_{\sx\sy}|_{x=y}\eer}tang_j[\ltSi_j]
\eeq
where  $tang_j$ is the degree of tangency of $\kill T_\sy$ and $\ltSi_{\sx\sy}$ along $\ltSi_j$.
Let $q=min(tang_j)$. In many cases  $tang_j>q$ only for a very few components.
Then we can write:
\beq
[\ltSi_{\sx\sy}]\Big([\kill T_\sy]-q[E_{x=y}]\Big)=[\ltSi_{\sx\sy'}]+
\sum_{\ber components~of\\\ltSi_{\sx\sy}|_{x=y}\eer}(tang_j-q)[\ltSi_j]
\eeq
where $E_{x=y}=\ltSi_{\sx\sy}|_{x=y}$ is the "exceptional divisor". This form simplifies the computations.
\subsection{On the validity range of results}\label{Sec Method Validitiy of Results}
From the explanation of the method it is
clear that an (obvious) necessary condition for the method to be applicable is: the stratum $\Si_{\sx\sy}$ is
smooth, irreducible of expected dimension.

A cheap sufficient condition is: $d\ge o.d.(\sx)+o.d.(\sy)$ (here $o.d.(\mS)$ is the degree of determinacy,
cf. the end of \S\ref{Sec Background SingularityTypes}). The bound can be obtained as follows.
Start from the stratum $\ltSi_{\sx\sy}$ and degenerate to the types: \omp s of multiplicities $o.d(\sx),o.d.(\sy)$.
 The formulas
certainly hold for $\ltSi_{\sx\sy}$ if the corresponding formulas hold for the degenerated stratum.
And in this later case an easy sufficient condition is
 $d\ge o.d.(\sx)+o.d.(\sy)$ (cf. theorem \ref{ThmTwoOrdinaryMultiplePoints}).

Of course, this sufficient bound is far from being necessary.
\section{Examples}\label{Sec Examples}
\subsection{Two \omp s}\label{Sec Example Two OMP's}
Let $\sx=x^{p+1}_1+x^{p+1}_2$, $\sy=y^{q+1}_1+y^{q+1}_2$, $p\ge q$.
As in the example \ref{ExOneOrdinaryMultiplePoint}
the natural candidate for the lifting of $\lSi_{\sx\sy}$ is the variety of triples $(x,y,f)$ with
$f$ having $\sy$ at $y$ and $\sx$ at $x$. To simplify the embedding we
blowup the auxiliary space $\mP^2_x\times\mP^2_y$ along the diagonal $\De=\{x=y\}$.
Geometrically we add the line $l=\overline{xy}$ (defined by a one-form, a point in the dual plane $\tP^2_l$).
The exceptional divisor is $E=\{(x,y,l)|~x=y,~l(x)=0\}$.
Thus the lifted stratum is defined as the strict transform:
\beq\ber\label{EqTwoMultiplePointsLiftedStratum}
\ltSi_{\sx\sy}(x,y,l):=\overline{\Big\{(x,y,l,f),~~x\ne y|~~f|_x^{(p)}=0=f|_y^{(q)},~~~l(x)=0=l(y)\Big\}}
\subset \Aux\times|\cO_{\mP^2}(d)|
\\
\Aux=\{(x,y,l)|l(x)=0=l(y)\}\subset\mP^2_x\times\mP^2_y\times\tP^2_l
\eer\eeq
\subsubsection{Stepwise intersection}\label{Sec Examples 2 OMP stepwise intersection}
\bthe
The projection $\ltSi_{\sx\sy}(x,y,l)\ra \Aux\subset \mP^2_x\times\mP^2_y\times\tP^2_l$
is a projective fibration over a smooth base, the projectivization of a vector bundle.
In particular the lifted stratum is a smooth locally complete intersection.
\ethe
\bpr
As the degree is high, the defining conditions are transversal outside the diagonal $x=y$.
(To check a specific fibre, fix the points $x,y$ then the conditions are just linearly independent
equations in $|\cO_{\mP^2}(d)|$.)
 Therefore the topological closure $\ltSi_{\sx\sy}(x,y,l)$ is an
irreducible reduced algebraic variety.

For $x\neq y$ the fibre of $\ltSi_{\sx\sy}(x,y,l)\ra \Aux$ over $(x,y,l)$ is a linear subspace of $|\cO_{\mP^2}(d)|$.
Consider the fibres over the diagonal. The points of $\ltSi_{\sx\sy}(x,y,l)|_{x=y}$
 correspond to collision of the two \omp s.
 As is shown in \S\ref{Sec Background Collisions With an OMP}, for $\sx\sy$ \omp s, the restriction is irreducible:
$\ltSi_{\sx\sy}(x,y,l)|_{x=y}=\ltSi_{\mS}(x,l,f)$. Here $\mS$ is the singularity type of
$(x^{p-q}_1+x^{p-q}_2)(x^{q+1}_1+x^{2q+2}_2)$
 and $l$ is the line tangent to $(x^{q+1}_1+x^{2q+2}_2)$.
Note that the fibre over $(x=y,l)$ is a linear subspace of $|\cO_{\mP^2}(d)|$.
Its co-dimension (computed e.g. as the number of $\mZ^2_{\ge0}$ points strictly below the \ND)
is $\bin{p+2}{2}+\bin{q+2}{2}$, i.e. precisely the codimension of the general fibre over $(x,y,l)$.
\\
Hence, as in the proof of proposition \ref{Thm Desingularization fo Strata}, there is the
natural morphism: $\ltSi_{\sx\sy}\to \Aux\to Gr(|\mP^n,\cO_{\mP^2}(d)|)$, giving to $\ltSi_{\sx\sy}$
the structure of projective fibration as on the diagram.
\beq
\bM \ltSi_{\sx\sy}(x,y,l)&\into&\tau&\sset Gr(\mP^n,|\cO_{\mP^2}(d)|)\times|\cO_{\mP^2}(d)|
\\\da &&\da&
\\\Aux&\into&Gr(\mP^n,|\cO_{\mP^2}(d)|)\eM
\eeq
Here $\tau$ is the tautological fibration,
the fibre over $[\mP^n]\in Gr$ is $\mP^n\sset|\cO_{\mP^2}(d)|$.
\epr
\bthe\label{ThmTwoOrdinaryMultiplePoints}
For $p\ge q$ the cohomology class of the lifted stratum $[\ltSi_{\sx\sy}(x,y,l)]\in
H^*(\mP^2_x\times\mP^2_y\times\tP^2_l\times|\cO_{\mP^2}(d)|,\mZ)$ is given by:
\beq
[\ltSi_{\sx\sy}(x,y,l)]=(L+X)(L+Y)\Big(F+(d-p)X\Big)^{p+2\choose{2}}\prod^q_{i=0}\prod^{q-i}_{j=0}
\Big(F+(d-i-j)Y+iX-jL-(p+1+i-j)E\Big)
\eeq
here $E=X+Y-L$ is the class of the exceptional divisor. (For the notations of the cohomology generators
cf. \S\ref{Sec Background Coordinates}). The formula is applicable for $d\ge p+q+2$.
\ethe
To get the solution of the enumerative problem (i.e. the degree of $\lSi_{\sx\sy}$) we should apply the Gysin homomorphism
corresponding to the projection $\ltSi(x,y,l)\ra\lSi$. In this case it means just to extract the coefficient
of $X^2Y^2L^2$. In
the case $\sx=\sy$ (i.e. $p=q$) the resulting answer should be also divided by 2, as the singular points are not ordered.
\bcor\label{ThmTwoOrdinaryMultPointNumericalCorollary}
In a few simplest cases the degree of $\lSi_{\sx\sy}$ is:
\li $\ber q=1.\\
deg(\lSi_{x^{p+1}_1+x^{p+1}_2,A_1})=\eer 9{p+3\choose{4}}(d-p)^3(d+p-2)-\frac{3}{4}{p+2\choose{3}}(10p^2+39p+7)(d-p)^2+
3{p+2\choose{3}}(d-p)(6+5p)
$
\\\\\li $\ber q=2.\\
\deg(\lSi_{x^{p+1}_1+x^{p+1}_2,D_4})=\eer\ber 45{p+3\choose{4}}(d-p)^3(d+p-4)+2(d-p)(8+3p+p^2)(35p^2+20p-12)\\-
\frac{5}{8}(p+1)(d-p)^2(14p^4+105p^3+147p^2+114p-80)
-6(85p^2+45p-28)\eer
$
\\\\\li $\ber q=3.\\
\deg(\lSi_{x^{p+1}_1+x^{p+1}_2,X_9})=\eer\ber 135{p+3\choose{4}}(d-p)^3(d+p-6)+2(d-p)(16+3p+p^2)(270p^2-20p-117)\\
-\frac{5}{8}(d-p)^2(54p^5+527p^4+948p^3+1853p^2-894p-1152)-14(830p^2-105p-348)\eer
$
\ecor
For $p=q=1$ this gives the classical result \cite[pg. 2]{Alufi98}.
For $p=2,q=1$ this coincides the results of \cite{KleiPien98} and \cite{Kaz03-hab}. For other cases the results seem to be new.
\beR\label{Rem Class depends non algebraically}
It is immediate from (\ref{EqTwoMultiplePointsLiftedStratum}) that the stratum is symmetric with respect
to the permutation $(x,p)\leftrightarrow(y,q)$. So, its
cohomology class is symmetric in $(X,p)\leftrightarrow(Y,q)$. The answer in \ref{ThmTwoOrdinaryMultiplePoints}
is obtained (in a very non-symmetric way) under the assumption $p\ge q$.
One might hope that this answer is still symmetric. However this happens for $(p,p)$ and $(p+1,p)$ cases only.
So, for the general answer one should substitute $max(p,q)$ and $min(p,q)$, i.e. the degree
is written as $\deg\lSi_{max(\sx,\sy),min(\sx,\sy)}$ with the obvious order on the types of \omp s.

This contradicts the naive expectation of algebraicity (as mentioned in \S\ref{Sec Intro Our Results}):
in general the cohomology classes of equisingular strata depend non-algebraically on the classical
singularity invariants.
\eeR
{\bf proof of the theorem:}\\
As in \S\ref{Sec Method Stepwise Intersection} represent $\ltSi_{\sx\sy}(x,y,l)$ as a stepwise intersection of the stratum $\ltSi^{(y)}_{\sx}(x,y,l)$ with hypersurfaces
 comprising the conditions $f|_y^{(q)}=0$. Start from the case of one singular point
\beq
\ltSi^{(y)}_{\sx}(x,y,l):=\Big\{(x,y,l,f)|~~f|_x^{(p)}=0,~~~l(x)=0=l(y)\Big\}\subset |\cO_{\mP^2}(d)|\times \Aux
\eeq
At j'th step  we have a variety $M_j\into|\cO_{\mP^2}(d)|\times \Aux$ and a hypersurface $V_{j+1}\subset|\cO_{\mP^2}(d)|\times \Aux$. Think about
the intersection $M_j\cap V_{j+1}$ as the pullback $i^*(V_{j+1})$. Then the task is to take the {\it
strict transform} of $V_{j+1}$, i.e. to subtract from the pullback the part of the
exceptional divisor over $\{x=y\}$.

The straightforward approach is just to consider the components of the tensor $f|_y^{(q)}$ (i.e. all the partials) and
intersect $\ltSi^{(y)}_{\sx}(x,y,l)$ with the corresponding hypersurfaces:  $\{\di^{n_0}_0\di^{n_1}_1\di^{n_2}_2 f|_y=0\}_{n_0+\dots n_2=q}$.
This will bring various complicated residual pieces. Instead, we represent these conditions as follows (cf. \S\ref{Sec Background Coordinates}
for notations):
\beq
\Big\{f|_y^{(i)}(\underbrace{x\dots x}_{i})=0\Big\}^q_{i=0},~~~~
\Big\{f|_y^{(i+1)}(\underbrace{x\dots x}_{i}\tv)=0\Big\}^{q-1}_{i=0},~~\dots .~~
\Big\{f|_y^{(i+j)}(\underbrace{x\dots x}_{i}\underbrace{\tv\dots \tv}_{j})=0\Big\}^{q-j}_{i=0},~~~~
f|_y^{(q)}(\underbrace{\tv\dots \tv}_{q})=0
\eeq
Here $\tv$ is a fixed point, so that the points $x,y,\tv$ generically do not lie on one line.
By direct check it is verified that for generic parameters (i.e. $y\ne x,~~\tv\notin \overline{xy}=l$)
these conditions are equivalent to $f|_y^{(q)}=0$. For non-generic situation each such equation will give a
reducible hypersurface in $M_j$, correspondingly a residual term should be subtracted.

\li $M_0:=\overline{\ltSi_{\sx}(x,y,l)\underset{x\ne y}{\cap}\{f|_y=0\}}$. The pullback of the hypersurface $\{f|_y=0\}$
to $\ltSi_{\sx}(x,y,l)$ consists of
the strict transform (the closure of the part over $x\ne y$) and the exceptional divisor $E$ (over $x=y$).
To calculate the multiplicity, expand $y=x+\ep v$, correspondingly:
\beq\label{Eq Proof Two OMP calc multipl}
0=f|_y=\underbrace{f|_x+\dots +\ep^pf|_x^{(p)}(v\dots v)}_{vanish}+\ep^{p+1}f|_x^{(p+1)}(\underbrace{v\dots v}_{p+1})+\dots
\eeq
i.e. the exceptional divisor enters with the multiplicity $(p+1)$. So, the strict transform is
$(f|_y=0)-(p+1)E$ and the total cohomology class:
$[\ltSi_{\sx}(x,y,l)]\Big([f|_y=0]-(p+1)[E]\Big)\in H^*(|\cO_{\mP^2}(d)|\times \Aux)$.

The points of $M_0$ satisfy:
\ls for $y\ne x$: $f|_x^{(p)}=0$ and $f|_y=0$
\ls for $y=x$: $f|_x^{(p)}=0=f|_x^{(p+1)}\underbrace{(v\dots v)}_{p+1}$
\li In the same way do all the intersections with $\Big\{f|_y^{(i)}(\underbrace{x\dots x}_{i})=0\Big\}^q_{i=1}$
(i.e. $M_i:=\overline{M_{i-1}\underset{x\ne y}{\cap}\big\{f|_y^{(i)}(\underbrace{x\dots x}_{i})=0\big\}}$).
At each step subtract the exceptional divisor with the necessary multiplicity. The resulting cohomology class is:
\beq
[\ltSi_{\sx}(x,y,l)]\prod^q_{i=0}\Big([f|_y^{(i)}(\underbrace{x\dots x}_{i})=0]-[(p+1+i)E]\Big)
\eeq
The points of $M_q$ satisfy:
\ls for $y\ne x$: $f|_x^{(p)}=0$ and $f|_y=0=f|_y^{(1)}(x)=\dots f|_y^{(q)}\underbrace{(x\dots x)}_{q}$
\ls for $y=x$: $f|_x^{(p)}=0=f|_x^{(p+1)}\underbrace{(v\dots v)}_{p+1}=\dots =f|_x^{(p+q+1)}\underbrace{(v\dots v)}_{p+q+1}$
\li Intersect with $\Big\{f|_y^{(i+1)}(\underbrace{x\dots x}_{i}\tv)=0\Big\}^{q-1}_{i=0}$.
Here $\tv$ is a fixed point in $\mP^2$. The intersection of each such hypersurface with the previously
obtained variety (say $M_j$) is reducible: it contains a component over $\tv\in\lxy=l$. (Note
that this component is a divisor.)
So, one should
take the "strict transform": $M_{j+1}=\Big(\{f|_y^{(i+1)}(\underbrace{x\dots x}_{i}\tv)=0\}-\{\tv\in\lxy\}\Big)\cap M_j$.
Now, every point of $M_{j+1}$ satisfies: $f|_y^{(i+1)}(\underbrace{x\dots x}_{i})=0$.
\\
Now, check the situation over the diagonal $x=y$, subtract its contribution as in the previous cases:
\beq
[\ltSi_{\sx}(x,y,l)]\prod^q_{i=0}\Big([f|_y^{(i)}(\underbrace{x\dots x}_{i})=0]-[(p+1+i)E]\Big)
\prod^{q-1}_{i=0}\Big([f|_y^{(i+1)}(\underbrace{x\dots x}_{i}\tv)=0]-[\tv\in l]-[(p+i)E]\Big)
\eeq
Here $\tv\in\overline{xy}=l$ is a divisor (all the lines in the plane through the fixed point $\tv$). Its cohomology class: $[\tv\in l]=L$.
It should be subtracted with multiplicity 1. To check this expand as previously (eq. \ref{Eq Proof Two OMP calc multipl}):
$\tv=v+\ep \tv_1$ (with $v\in l$, $\tv_1\notin l$), then
\beq
f|_y^{(i+1)}(\underbrace{x\dots x}_{i}\tv)=\underbrace{f|_y^{(i+1)}(\underbrace{x\dots x}_{i}v)}_{vanishes}+\ep f|_y^{(i+1)}(\underbrace{x\dots x}_{i}\tv_1)
\eeq

The points of the so obtained variety satisfy:
\ls for $y\ne x$: $f|_x^{(p)}=0$ and $\{f|_y^{(i)}\underbrace{(x\dots x)}_{i-1}=0\}^q_{i=1}$
\ls for $y=x$: $f|_x^{(p)}=0=f|_x^{(p+1)}\underbrace{(v\dots v)}_{p}=\dots =f|_x^{(p+q)}\underbrace{(v\dots v)}_{p+q-1}=f|_x^{(p+q+1)}\underbrace{(v\dots v)}_{p+q+1}$
\li Do the rest of intersections, at each step subtracting (with appropriate multiplicities) the
exceptional divisor and the class $[\tv\in l]$. Finally we get:
\beq
[\ltSi_{\sx\sy}(x,y,l)]=[\ltSi_{\sx}(x,y,l)]\prod^q_{i=0}\prod^{q-i}_{j=0}
\Big([f|_y^{(i+j)}(\underbrace{x\dots x}_{i}\underbrace{\tv\dots \tv}_{j})=0]-j[\tv\in l]-[(p+1+i-j)E]\Big)
\eeq
The points of this variety satisfy:
\ls for $y\ne x$: $f|_x^{(p)}=0=f|_y^{(q)}$
\ls for $y=x$: $f|_x^{(p)}=0=f|_x^{(p+1)}\underbrace{(v\dots v)}_{p+1-q}=f|_x^{(p+2)}\underbrace{(v\dots v)}_{p+3-q}=\dots =f|_x^{(p+q+1)}\underbrace{(v\dots v)}_{p+q+1}$

Substitute now the cohomology classes for the conditions (cf. \S\ref{Sec Background}).
(Note that as $\tv$ is a fixed point, the condition $\tv\in l$ is just one linear condition on $l$,
its class is $L$. The class of the exceptional divisor is $[E]=X+Y-L$, cf. \S\ref{Sec Background Coordinates}.)
This proves the formula.
\\\\
Finally, regarding the sufficient bound for validity of the formula. Note that all the conditions that have appeared in
the proof concern only the $(p+q+1)$'th jet of the function $f$. So, the bound $d\ge p+q+2$ is sufficient.
\epr
\subsubsection{Degenerations: chipping off a line}\label{Sec Intro Example 2OMP Degen to Reducible}~\\
\parbox{14cm}
{Use the following elementary observation. Let $C^{(d)}_{x_{p+1},y_{q+1}}$ be a curve of degree $d$ with
the points $x_{p+1},y_{q+1}$ of multiplicities $p+1,q+1$. Let $l=\overline{x_{p+1},y_{q+1}}$.
Suppose the intersection multiplicity satisfies: $|l\cap C^{(d)}_{x_{p+1},y_{q+1}}|>d$.
Then the curve is reducible: $C^{(d)}_{x_{p+1},y_{q+1}}=l\cup C^{(d-1)}_{x_{p},y_{q}}$.

We degenerate to get to this situation. The degeneration is done by forcing one of the branches at $x$ to be tangent to $l$ with
high degree of tangency, i.e. "killing" the corresponding monomials (cf. the \ND).
}
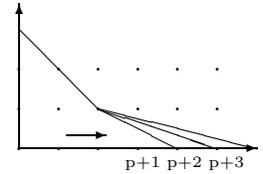
\begin{picture}(0,0)(-10,20)
\mesh{0}{0}{5}{2}{15}{90}{55}
\put(0,45){\line(1,-1){30}}\put(30,15){\line(2,-1){30}}\put(30,15){\line(3,-1){45}}\put(30,15){\line(4,-1){60}}
\put(18,5){\vector(1,0){15}}
\put(40,-7){\tiny{p+1~~p+2~~~p+3}}
\end{picture}
\\
This amounts to intersection of the stratum with the corresponding "degenerating hypersurfaces". Their cohomology
classes are given in \S\ref{Sec Method CohomologyClassDegenerations}.
Let $\ltSi^{(d)}_{p+1,q+1}(x,y,l)$ denote the lifted stratum of degree $d$ curves with two \omp s,
of multiplicities $i,j$.
\bthe The intersection of
the stratum $\ltSi^{(d)}_{p+1,q+1}(x,y,l)$ with $(d-q-p-1)$ degenerating hypersurfaces results in the stratum of reducible curves
$\Xi^{(d-1)}_{p,q}:=\Big\{(x,y,l,C_d)\Big|~C_d=l\cup C_{d-1},~~(x,y,l,C_{d-1})\in\ltSi^{(d-1)}_{p,q}(x,y,l)\Big\}$.
In particular, the cohomology classes of the two strata are related by the equation:
\beq
[\ltSi^{(d)}_{p+1,q+1}(x,y,l)]\prod^{d-q-1}_{k=p+1}\Big(F+(d-k)X+kY-(k+q+1)E\Big)=[\Xi^{(d-1)}_{p,q}]\in
 H^*(\mP^2_x\times\mP^2_y\times\tP^2_l\times|\cO_{\mP^2}(d)|),~~~here~~E=X+Y-L
\eeq
\ethe
\bpr
Start from the stratum $\ltSi^{(d)}_{p+1,q+1}(x,y,l)$, the corresponding curves satisfy:
$f^{(p)}|_x=0=f^{(q)}|_y$. Consider the degenerating hypersurfaces:
\beq
D_{p+1}:=\{f|_x^{(p+1)}\underbrace{(y\dots y)}_{p+1}=0\},~~D_{p+2}:=\{f|_x^{(p+2)}\underbrace{(y\dots y)}_{p+2}=0\},\dots .
D_{d-q-1}:=\{f|_x^{(d-q-1)}\underbrace{(y\dots y)}_{d-q-1}=0\}
\eeq
These conditions, imposed onto the curve $C^{(d)}_{x_{p+1},y_{q+1}}$, force the intersection multiplicity:
$mult_x(C^{(d)}_{x_{p+1},y_{q+1}},\overline{xy})\ge (p+1)+(d-p-q-1)=d-q$.
(Pass to the local coordinate system to verify that these are precisely the derivatives in the given direction.)

Hence, outside the diagonal $\{x=y\}$ one has:
\beq
\ltSi^{(d)}_{p+1,q+1}(x,y,l)\underset{x\neq y}{\cap} \bigcap_{j=p+1}^{d-q-1}D_j=\Xi^{(d-1)}_{p,q}\cap \{x\neq y\}
\eeq
The cohomology classes are (see \S\ref{Sec Method CohomologyClassDegenerations}):
$[D_j]=F+(d-j)X+jL\in H^2(\mP^2_x\times\mP^2_y\times\tP^2_l\times|\cO_{\mP^2}(d)|)$.

Consider the intersection near the diagonal. Recall from\S\ref{Sec Background Collisions}, that the points over the
diagonal $\ltSi^{(d)}_{p+1,q+1}(x,y,l)|_{x=y}$ correspond to the type with the representative:
$(x^{p-q}+y^{p-q})(x^{q+1}+y^{2q+2})=x^{p+1}+x^{q+1}y^{p-q}+y^{p+q+2}$.
The defining equations for such a fibre are:
\beq
f|_y^{(p)}=0,~~f|_y^{(p+1)}(\underbrace{v\dots v}_{p-q+1})=0,~~\dots ~~ f|_y^{(p+q+1)}(\underbrace{v\dots v}_{p+q+1})=0
\eeq
here $v$ is a point on the line $l$, $x\neq v$. Thus the intersection $D_{p+1}\cap\ltSi^{(d)}_{p+1,q+1}(x,y,l)$
contains $\ltSi^{(d)}_{p+1,q+1}(x,y,l)|_{x=y}$ as an irreducible component. To find the multiplicity of this
component, expand $x=y+\ep v$:
\beq
f|_x^{(p+1)}\underbrace{(y\dots y)}_{p+1}=\underbrace{\dots .+
\ep^{q}f|_y^{(p+q+1)}(\underbrace{v\dots v}_{p+q+1})}_{vanish~on~\ltSi^{(d)}_{p+1,q+1}(x,y,l)|_{x=y}}+
\ep^{q+1}f|_y^{(p+q+2)}(\underbrace{v\dots v}_{p+q+2})+\dots
\eeq
Therefore for the cohomology classes have:
\beq\ber
[\ltSi^{(d)}_{p+1,q+1}(x,y,l)][D_{p+1}]=
[\overline{\lSi^{(d)}_{p+1,q+1}(x,y,l)\underset{x\neq y}{\cap}D_{p+1}}]+(p+q+2)[\ltSi^{(d)}_{p+1,q+1}(x,y,l)|_{x=y}]
\\
i.e. [\ltSi^{(d)}_{p+1,q+1}(x,y,l)]\Big([D_{p+1}]-(p+q+2)E\Big)=
[\overline{\lSi^{(d)}_{p+1,q+1}(x,y,l)\underset{x\neq y}{\cap}D_{p+1}}]
\eer\eeq
which is the first term in the product of the statement. Continue by induction.

After the k'th step the points of the stratum
$\overline{\lSi^{(d)}_{p+1,q+1}(x,y,l)\underset{x\neq y}{\cap}\capl^k_{i=1}D_{p+i}}$
satisfy:
\beq
f^{(p)}|_x=0=f^{(q)}|_y,~~~f|_x^{(p+1)}\underbrace{(y\dots y)}_{p+1}=0, ~ \dots  ~ f|_x^{(p+k)}\underbrace{(y\dots y)}_{p+k}=0
\eeq
To check the intersection of $D_{p+k+1}$ on the diagonal we should obtain the defining equations of
 $\overline{\lSi^{(d)}_{p+1,q+1}(x,y,l)\underset{x\neq y}{\cap}\capl^k_{i=1}D_{p+i}}$ near $x=y$. This is done
 by the flat limit. The procedure similar to that of \ref{Sec Background Collisions} gives
 the degree of the intersection $D_{p+k+1}$ with the diagonal: $p+q+2+k$. So, the k'th step
 brings to the product the factor: $[D_{p+k+1}]-(p+q+2+k)E$.
\epr
\\\\
By applying this process $q+1$ times one arrives at a curve with an \omp~ of multiplicity $(p+1)$
and a line of multiplicity $q+1$ as a component.
For such a stratum the answer is known classically. As every intersection is invertible this solves
the enumerative problem.

\subsection{The class of $\ltSi_{\sx,A_1}$ by degenerations to \omp}\label{Sec Example Degeneration to OMP's}
We consider here the simplest case: $\sy=A_1$, then the degeneration process of $\ltSi_{\sx,A_1}$
is quite explicit.

Let $T_\sx=l_1^{p_1}\dots l_k^{p_k}$ be the tangent cone of $\sx$ and $\ltSi_\sx(x,l_1\dots l_k)$
the corresponding lifting, as in \S\ref{Sec Method Strata of Unisingular Curves}. So, we start from the lifted stratum:
\beq
\ltSi_{\sx,A_1}(x,y,l,l_1\dots l_k):=\overline{\Big\{\ber (x,y,l,l_1\dots l_k,C)\\ x\neq y, ~ ~ l=\overline{xy}\eer\Big|
\ber \text{ $C$ has $\sx$ at $x$ and $A_1$ at $y$} \\ T_{(C,x)}=l_1^{p_1}\dots l_k^{p_k}\eer\Big\}}
\eeq
Intersect with the divisor $\kill T_\sx$, cf.\S\ref{Sec Method Killing Tangent Cone},
such that $\ltSi_\sx(x,l_1\dots l_k)\cap \kill T_\sx=\ltSi_{\sx'}(x,l_1\dots l_k)$. Here $\sx'$ is
the degenerated singularity type, recall that if $\sx$ is linear singularity then $\sx'$ is well defined.

Note that $T_\sx\neq T_{\sx'}$, hence some of $l_i$ of $\ltSi_{\sx'}(x,l_1\dots l_k)$ can be
not related to the parameters of the singularity $\sx'$. In the extremal case, when $\sx'$ is an \omp, one has:
$\ltSi_{\sx'}(x,l_1\dots l_k)=\ltSi_{\sx'}(x)\cap\bigcap_i\{x\in l_i\}$.

The intersection $\ltSi_{\sx,A_1}(x,y,l,l_1\dots l_k)\cap \kill T_\sx$ brings in addition the residual pieces over
the diagonal:
\bthe Killing the tangent cone results in the cohomological equation in $H^*(\Aux\times|\cO_{\mP^2}(d)|,\mZ)$
\beq
\Big[\ltSi_{\sx,A_1}(x,y,l,\{l_i\})\Big][\kill T_\sx]=
\Big[\ltSi_{\sx',A_1}(x,y,l,\{l_i\})\Big]+2[x=y]\Big[\ltSi_{\mS_1}(x,l,\{l_i\})\Big]+
[x=y]\suml_{\tinyA\bM i\text{ such that}\\p_i=1\eM}[l=l_i]\Big[\ltSi_{\mS_2}(x,l,\{l_i\})\Big]
\eeq
where the strata for $\mS_1,\mS_2$ are defined by
\beq\ber
\ltSi_{\mS_1}(x,l,\{l_i\})=
\ltSi_{\mS'}(x,\{l_i\})\cap\Big\{x\in l, ~ ~ f|_x^{(p)}=0=f|_x^{(p+1)}(\underbrace{v\dots v}_{p})~\forall v\in l\Big\}
\\
\ltSi_{\mS_2}(x,l,\{l_i\})=
\ltSi_{\mS'}(x,\{l_i\})\cap\Big\{x\in l, ~ ~ f|_x^{(p)}=0=f|_x^{(p+1)}(\underbrace{v\dots v}_{p+1})~\forall v\in l\Big\}
\eer\eeq
The needed class $[\ltSi_{\sx,A_1}(x,y,l,\{l_i\})]$ is determined from this equation  uniquely.
\ethe
This theorem expresses the class $\Big[\ltSi_{\sx,A_1}(x,y,l,\{l_i\})\Big]$ in terms of
the degenerated stratum $\Big[\ltSi_{\sx',A_1}(x,y,l,\{l_i\})\Big]$ and the known classes
 $\Big[\ltSi_{\mS_1}(x,l,\{l_i\})\Big]$ and $\Big[\ltSi_{\mS_2}(x,l,\{l_i\})\Big]$, the later
  two correspond to curves with one singular point (of a linear type). Hence, after a few degenerating
  steps (not more than the order of determinacy of $\sx$ minus the multiplicity of $\sx$),
 one arrives at the stratum of \omp s, the class of which was computed in the preceding section.

In particular, if $o.d.(\sx)-\mult(\sx)=1$, i.e. $\sx$ has a representative of the form $f_p+f_{p+1}$,
with $f_{p+1}$ generic, then the problem is solved in one step.
If $\sx$ is already an \omp~ then the theorem can be used to compute the class
$[\ltSi_{\sx',A_1}(x,y,l)]$ in terms of $[\ltSi_{\sx,A_1}(x,y,l)]$, i.e. it gives an alternative (and fast)
way to solve the enumerative problem for \omp s.

Various numerical answers are given in Appendix.
\\
\bpr We need only to understand the residual pieces over the diagonal $\{x=y\}$ and their multiplicities.
Consider the locally defining equations of $\ltSi_{\sx A_1}(x,y,l,\{l_i\})$
off $x=y$. As in \S\ref{Sec Background Collisions} take the limit $yto x$, i.e.
expand $y=x+\ep v$. Then the preliminary set of equations
consists of the equations of $\sx$ and local expansion of $\sy=A_1$:
\beq\ber
f^{(p)}|_x\sim SYM(l^{p_1}_1,\dots ,l^{p_k}_k), \dots .~ \text{ further equations of } \sx, ~\dots .
\\
f^{(p)}|_x(\underbrace{v\dots v}_{p-1})+\ep f^{(p+1)}|_x(\underbrace{v\dots v}_{p})+\dots =0=
f^{(p+1)}|_x(\underbrace{v\dots v}_{p+1})+\ep f^{(p+2)}|_x(\underbrace{v\dots v}_{p+2})+\dots
\eer\eeq
Here the proportionality of tensors is the only condition of $\sx$ involving $f^{(p)}|_x$.

This equation gives rise to the syzygies. We should take the flat limit as $\ep\to0$.
The following cases are possible:
\li $\exists i:$ $p_i\ge2$ and $l_i=l$. Then in the equation above the term $f^{(p)}|_x(\underbrace{v\dots v}_{p-1})$
necessarily vanishes (as $v\in l_i$). So, the equations place no conditions on $f^{(p)}|_x$ except
for $f^{(p)}|_x\sim SYM(l^{p_1}_1,\dots ,l^{p_k}_k)$. Thus, generically $f^{(p)}|_x\neq0$ and this locus of $\Aux$
does not contribute to the cohomology class (cf.\S\ref{Sec Residual Pieces sx to OMP}).

\li $\exists i:$ $p_i=1$ and $l_i=l$, but $\{l_i\}$ are pairwise distinct. Then in the equations above the
number of independent (linear) conditions on $f^{(p)}|_x$ equals precisely the number
of independent entries. So, the equations of this piece are $f^{(p)}|_x=0=f^{(p+1)}|_x(\underbrace{v\dots v}_{p+1})$
and some equations of $\sx$ on higher derivatives.

\li $\forall i$: $l_i\neq l$ and $\{l_i\}$ are distinct. Then the number of independent (linear) conditions on
$f^{(p)}|_x$ is bigger by one than the number of its independent entries. So,
 the equations of the contributing residual piece are:
$f^{(p)}|_x=0=f^{(p+1)}|_x(\underbrace{v\dots v}_{p})$.

Note that in both cases the residual pieces are reduced, as all the equations are linear in $f$ or its derivatives.
\\
\\
Now compute the tangency degrees of the divisor $\kill T_\sx$ and the residual pieces.
Namely, add to the locally defining ideal of $\ltSi_{\sx A_1}(x,y,l,\{l_i\})$ near the
contributing residual piece the local equation of $\kill T_\sx$ and check the primary decomposition of
the bigger ideal. In the two relevant cases we have:
\li $\exists i:$ $p_i=1$ and $l_i=l$, but $\{l_i\}$ are pairwise distinct. Then the enlarged ideal is generated by
\beq
f^{(p)}|_x=0, ~ ~\ep f^{(p+1)}|_x(\underbrace{v\dots v}_{p})+\ep^2 f^{(p+2)}|_x(\underbrace{v\dots v}_{p+1})\dots , ~ ~
f^{(p+1)}|_x(\underbrace{v\dots v}_{p+1})+\ep f^{(p+2)}|_x(\underbrace{v\dots v}_{p+2})+\dots ,
\dots \text{ conditions on higher derivatives}\dots
\eeq
And its primary decomposition consists of the two parts:
$\Big\{f^{(p)}|_x,$ $\ep,$ $f^{(p+1)}|_x(\underbrace{v\dots v}_{p+1}),\dots \Big\}$ and
$\Big\{f^{(p)}|_x,$ $f^{(p+1)}|_x(\underbrace{v\dots v}_{p})+\ep^2 f^{(p+2)}|_x(\underbrace{v\dots v}_{p+1})\dots ,\dots \Big\}$.
So the tangency in this case is 1.
\li $\forall i$: $l_i\neq l$ and $\{l_i\}$ are distinct. Then the enlarged ideal is generated by
\beq
f^{(p)}|_x=0, ~ ~
\ep^2 f^{(p+2)}|_x(\underbrace{v\dots v}_{p+1})+\dots =0=
f^{(p+1)}|_x(\underbrace{v\dots v}_{p})+\dots ,
\dots \text{ conditions on higher derivatives}\dots
\eeq
So the tangency in this case is 2.
\epr

\appendix
\section{Some numerical results for enumeration of curves with two singular points}\label{AppendixNumericalResults}
We give below explicit expressions $deg\lSi_{\sx,\sy}$ for various types $\sx$ and $\sy$ an \omp.
This is only a tiny amount of possible enumerative results, for more formulas
cf. the attached Mathematica file in the Arxiv \cite{Ker07Bisingular} or on author's homepage.

As is well known (cf. \cite{Kaz01,Kaz03-1}) for a collection of singularities $\mS_1\dots \mS_r$
the degree of the stratum $\lSi_{\mS_1\dots \mS_r}$ is expressed as:
\beq
\deg(\lSi_{\mS_1\dots \mS_r})=\frac{1}{|Aut|}\sum_{J_1\bigsqcup\dots \bigsqcup J_k}
S_{\{\mS_{i}\}_{i\in J_1}}\dots .S_{\{\mS_{i}\}_{i\in J_k}}
\eeq
Here the sum is over all the possible decompositions $\{\mS_1,\dots ,\mS_r\}=\bigsqcup J_i$
and $S_{\{\mS_{i}\}_{i\in J_1}}$ is the specialization of Thom polynomial, $S_{\sx}=\deg(\lSi_{\sx})$.
One divides by $|Aut|$ the cardinality of the automorphisms (depending on the coincidence of the types).
In particular, for two singular points: $\deg(\lSi_{\sx\sy})=\deg(\lSi_{\sx})\deg(\lSi_{\sy})+S_{\sx\sy}$.

The degrees $\deg(\lSi_{\mS_i})$ are known e.g. from \cite[\S A.2]{Ker06 enum.unising.curves},
for particular cases cf. example \ref{Sec Method Examples One Sing Point Many Branches}.
For completeness we present them here.

\li For the \omp: $\deg(\lSi_{x^{p+1}_1+x^{p+1}_2})=\bin{\bin{p+2}{2}}{2}(d-p)^2$.
\li For the singularity type $\sx$ with the representative $f_p+x^{p+1}_1+x^{p+1}_2$,
and the tangent cone $T_{\sx}=l^{p_1}_1\dots l^{p_k}_k$, with $\sum p_i=p$:
\beq\smallA
\deg(\lSi_{\sx})=\frac{1}{|G|}\Big(\prod p_i\Big)\Bigg(k!\bin{\bin{p+2}{2}-1-k}{2}(d-p)^2+
(k-1)!\bin{k+1}{2}\bin{\bin{p+2}{2}-2-k}{1}(d-p)\sum p_i
+(k-2)!\bin{k}{2}\bin{k+2}{2}\suml_{1\le i<j\le k}p_i p_j\Bigg)
\eeq
Here $G$ is the symmetry group, permuting lines with coinciding $p_j$'s.
\li For the singularity type $\sx$ with the representative $(x^{p-1}_1+x^p_2)(x_1+x^{2}_2)$,
\beq
\deg(\lSi_{\sx})=\frac{p(p+4)(p-1)}{8}(2p^3+7p^2-5p-2)(d-p)^2+(\bin{p+2}{2}-3)p^2(d-p)(d-2(p+1)),\\
\eeq
Below are the degrees of strata for curves with two singular points.
For two \omp s the final answers are also given in \S\ref{Sec Example Two OMP's}. Recall, for $\lSi_{\sx\sy}$
we always assume: $\mult(\sx)\ge \mult(\sy)$. If $\sx=\sy$ the final answer should be divided by 2.
\beq\ber
\deg(\lSi_{x^{p+1}_1+x^{p+1}_2,A_1})=\deg(\lSi_{x^{p+1}_1+x^{p+1}_2})\deg(\lSi_{A_1})-
\frac{3}{4}{p+2\choose{3}}(d-p)^2(3p+4)(p^2+3p+4)+3{p+2\choose{3}}(d-p)(5p+6)
\\
\deg(\lSi_{x^{p+1}_1+x^{p+1}_2,D_4})=\deg(\lSi_{x^{p+1}_1+x^{p+1}_2})\deg(\lSi_{D_4})+
\Bigg(\smallA\ber-\frac{5}{8}(d-p)^2(p+1)(3p-1)(p^2+3p+8)(p^2+3p+10)
\\+2(d-p)(p^2+3p+8)(35p^2+20p-12)-6(85p^2+45p-28)\Bigg)
\eer
\\
\deg(\lSi_{x^{p+1}_1+x^{p+1}_2,X_9})=\deg(\lSi_{x^{p+1}_1+x^{p+1}_2})\deg(\lSi_{X_9})
+\Bigg(\smallA\ber -\frac{5}{8}(d-p)^2(3p+2)(3p-2)(p^2+3p+16)(p^2+3p+18)\\
+2(d-p)(p^2+3p+16)(270p^2-20p-117)-14(830p^2-105p-348) \eer\Bigg)
\eer\eeq
\beq
\deg(\lSi_{x^{p}_1+x^{p+1}_2,A_1})=\deg(\lSi_{x^{p}_1+x^{p+1}_2})\deg(\lSi_{A_1})
-\frac{3}{8}p^4(3+p)(d-p)^2(p^2+3p-2)-\frac{3}{2}(p-1)p^3(d-p)(p^2+3p-2)+3p^4
\eeq
\beq
\deg(\lSi_{y(x^{p}_1+x^{p+1}_2),A_1})=\deg(\lSi_{y(x^{p}_1+x^{p+1}_2)})\deg(\lSi_{A_1})
+\Big(\smallA\ber-\frac{p^2(d-p)^2(5 + p)}{8}(2 + 5p + p^2)(2 + 11p + 6p^2)\\\frac{(p+1)(p+5)p^3}{4}(d-p)(2p+3)(3p+4)
-\frac{p^2(p-1)}{8}(6p^4+41p^3+55p^2+64p+92)\eer\Big)
\eeq
\beq
\deg(\lSi_{(x^2+y^2)(x^{p}_1+x^{p+1}_2),A_1})=\deg(\lSi_{(x^2+y^2)(x^{p}_1+x^{p+1}_2)})\deg(\lSi_{A_1})
+\Big(\smallA\ber -\frac{p(d-p)^2(1+p)(p+6)}{4}(p^2+7p+4)(9p^2+34p+24)+\\
+p(d-p)(p^2+7p+4)(9p^4+79p^3+220p^2+216p+63)\\
-p(9p^6+124p^5+587p^4+1316p^3+1480p^2+654p+60)\eer\Big)
\eeq
\beq
\deg(\lSi_{(x^{p-1}_1+x^p_2)(x_1-x^{2}_2),A_1})=\deg(\lSi_{(x^{p-1}_1+x^p_2)(x_1-x^{2}_2)})\deg(\lSi_{A_1})
+\Big(\smallA\ber -9{p+3\choose{4}}p(d-p)^2(4+p+2p^2)+
3p(p^4+3p^3+3p^2+4p-4)\\-\frac{3}{2}p^2(3+p)(d-p)(p^3-3p^2-p-8)\eer\Big)
\eeq
\beq
p>2, ~ \deg(\lSi_{(x^{p-2}_1+x^{p-2}_2)(x^2_1-x^{4}_2),A_1})=
\deg(\lSi_{(x^{p-2}_1+x^{p-2}_2)(x^2_1-x^{4}_2)})\deg(\lSi_{A_1})+
\Big(\tinyA\ber -\frac{3(d-p)^2}{8}(p^2+3p+6)(p^2+3p+8)(3p^3+6p^2+12p+5)\\
+\frac{3(d-p)}{2}(p^2+3p+6)(3p^4+23p^3+48p^2+81p+29)\\-3(57+182p+118p^2+51p^3+10p^4)\eer\Big)
\eeq
For small values $p=2,3$ these formulas reproduce the results of \cite{Kaz03-hab}.

\end{document}